\documentclass[11pt]{article}
\usepackage[T1]{fontenc}
\usepackage{latexsym,amssymb,amsmath,amsfonts,amsthm}
\usepackage{graphics}
\usepackage{graphicx}
\usepackage{mathrsfs}
\usepackage{subfigure}
\usepackage{color}
\usepackage{hyperref}
\newcommand{\R}{{\mathbb R}}
\newcommand{\Z}{{\mathbb Z}}
\newcommand{\N}{{\mathbb N}}
\newcommand{\C}{{\mathbb C}}
\newcommand{\Sp}{{\mathbb S}}
\newcommand{\s}{{\mathbb S}}

\newcommand{\no}{\nonumber}
\newcommand{\be}{\begin{eqnarray}}
\newcommand{\ben}{\begin{eqnarray*}}
\newcommand{\en}{\end{eqnarray}}
\newcommand{\enn}{\end{eqnarray*}}
\newcommand{\ba}{\backslash}
\newcommand{\pa}{\partial}

\newcommand{\ov}{\overline}

\newcommand{\G}{\Gamma}

\newcommand{\Om}{\Omega}

\newcommand{\ra}{\rightarrow}
\newtheorem{theorem}{Theorem}[section]

\newtheorem{remark}[theorem]{Remark}

\newtheorem{example}{Example}[]

\definecolor{xxx}{rgb}{0,0,0}%{0,0,0}¼´ÎªºÚÉ«
\definecolor{mgq}{rgb}{0,0.5,1}
\definecolor{hgh}{rgb}{1,0,0}

\begin{document}
\renewcommand{\theequation}{\arabic{section}.\arabic{equation}}
\begin{titlepage}

\title{Detection of a piecewise linear crack with one incident wave}

\author{Xiaoxu Xu\thanks{School of Mathematics and Statistics, Xi'an Jiaotong University, Xi'an, Shaanxi, 710049, China
(xuxiaoxu@xjtu.edu.cn).}
\and Guanqiu Ma\thanks{Corresponding author: School of Mathematical Sciences and LPMC, Nankai University, Tianjin
300071, China {(gqma@nankai.edu.cn)}.}
\and Guanghui Hu\thanks{School of Mathematical Sciences and LPMC, Nankai University, Tianjin
300071, China (ghhu@nankai.edu.cn).}}

\date{}
\end{titlepage}
\maketitle
\vspace{.2in}
\begin{abstract}
This paper is concerned with inverse crack scattering problems {\color{xxx}for  time-harmonic acoustic waves}.
We prove that a piecewise linear crack with the sound-soft boundary condition in {\color{xxx}two dimensions} can be uniquely determined by the far-field data corresponding to a single incident plane wave or point source.
We propose two non-iterative methods for imaging the location and shape of a crack.
The first one is a contrast sampling method, while the second one is a variant of the classical factorization method but only with one incoming wave.
Newton's iteration method is then employed for getting a more precise reconstruction result.
Numerical examples are presented to show the effectiveness of the proposed hybrid method.

\vspace{.2in} {\bf Keywords}: Helmholtz equation, crack, inverse scattering, single wave, uniqueness, hybrid method.
 \end{abstract}

\section{Introduction}

Inverse scattering problems aim to detect and identify the location, shape and physical properties of an unknown scatterer from the measurement data incited by incident waves.
The unknown target can be an impenetrable bounded obstacle, an inhomogeneous medium or an unbounded rough surface and so on.
The measurement data {\color{xxx}could be far-field data or near-field data.
The far-field data are far-field patterns which are also known as the scattering amplitudes, while the near-field data consist of the measurements of scattered waves or total waves}.
We refer the readers to \cite{Cakoni14,CK19} for an overview of inverse time-harmonic acoustic and electromagnetic scattering problems.
In this paper, we focus on inverse acoustic scattering by a sound-soft crack {\color{xxx}in $\R^2$}.
Let $\Gamma\subset\R^2$ be a piecewise linear crack embedded in an isotropic and homogeneous medium.
More precisely, $\G$ is supposed to be a piecewise linear curve lying on the boundary $\pa\Om$ of some convex polygon $\Om$.
Suppose the crack $\G$ is illuminated by the plane wave
\be\label{eq-inc}
u^i=u^i(x,d):=e^{ikx\cdot d},
\en
where $k\!>\!0$ is the wave number and $d\!\in\!\Sp$ is the incident direction with $\Sp\!:=\!\{x\in \mathbb{R}^2 \!:\!|x|\!=\!1\}$ denoting the unit circle.
Let $u^s(x,d)$ denote the scattered field. Then the total field $u(x,d)=u^i(x,d)+u^s(x,d)$ satisfies the equations:
\be\label{eq1}
\Delta u+k^2u=0 && \text{in}\;\R^2\ba\ov{\Gamma},\\ \label{eq2}
u_\pm=0 && \text{on}\;\Gamma,\\ \label{eq3}
\lim_{r\ra\infty}\sqrt{r}\left(\frac{\pa u^s}{\pa r}-iku^s\right)=0, && r=|x|,
\en
where (\ref{eq1}) is the so-called Helmholtz equation and (\ref{eq3}) is the Sommerfeld radiation condition that ensures the existence of a unique solution to (\ref{eq1})--(\ref{eq3}).
Since the crack is assumed to be sound-soft, the total field $u$ satisfies the homogeneous Dirichlet boundary condition (\ref{eq2}) on both sides of $\G$, where
\ben
{\color{xxx}u_\pm(x):=\lim_{t\ra+0}u(x\pm t\nu(x)),\quad x\in\G,}
\enn
with $\nu$ denoting the unit outward normal to $\pa\Om$.

The existence {\color{xxx}of a unique solution to the scattering problem (\ref{eq1})--(\ref{eq3}) has} been established in \cite[Section 8.7]{Cakoni14}.
By \cite[Theorem 2.6]{CK19}, the scattered field $u^s$ has the asymptotic behavior
\be\label{eq13}
u^s(x)=\frac{e^{ik|x|}}{\sqrt{|x|}}\left\{u^\infty(\hat x)+O\left(\frac1{|x|}\right)\right\},\quad|x|\ra\infty,
\en
uniformly for all observation directions $\hat x:=x/|x|$. Here, $u^\infty$ is called the far-field pattern of the scattered field $u^s$, which is an analytic function over the unit circle $\s$.
{\color{xxx}Given a point source}
\be\label{ps}
w^i(x,y):=\Phi_k(x,y),\quad y\in\R^2\ba\ov{\G},\quad x,y\in\R^2,
\en
we denote the scattered field, total field and its far-field pattern by $w^s(x,y)$, $w(x,y)$ and $w^\infty(\hat x,y)$, respectively.
Here, $\Phi_k(x,y)$ is the fundamental solution to the two-dimensional Helmholtz equation, i.e.,
\ben
\Delta_x\Phi_k(x,y)+k^2\Phi_k(x,y)=-\delta(x-y).
\enn
 It is well known that $\Phi_k(x,y)=\frac{i}4H_0^{(1)}(k|x-y|)$ with $H_0^{(1)}$ denoting the Hankel function of the first kind of order zero (see \cite[Section 3.5]{CK19}).

In this paper we are interested in the uniqueness and numerical method for inverse crack scattering problems with one incident wave.
More precisely, {\color{xxx}for a fixed wave number $k_0>0$} we want to reconstruct the location and shape of the crack $\G$ from {\color{xxx}a knowledge of the far-field pattern $\left\{u^\infty(\hat{x},d_0;k_0): \hat{x}\in\Sp\right\}$ for a fixed $d_0\in\Sp$ or} $\left\{w^\infty(\hat{x},y_0;k_0):\hat{x}\in\Sp\right\}$ for a fixed $y_0\in\R^2\ba\ov{\G}$. %Throughout this paper the wave number $k>0$ is {\color{xxx}arbitrarily} fixed.

Uniqueness in inverse scattering is concerned with the question whether the measurement data can uniquely determine the unknown target.
Assuming two different scatterers producing the same far-field patterns
for all incident directions, one can obtain a contradiction by considering a sequence of solutions with a singularity moving towards a boundary point of one scatterer that is not contained in the other scatterer (see \cite[Theorem 5.6]{CK19} in inverse obstacle scattering and \cite[Theorem 8.39]{Cakoni14} in inverse crack scattering).
If the scatterer is a convex polyhedral (or polygonal) obstacle of Dirichlet kind, the uniqueness to inverse acoustic and elastic scattering problems can be established with a single incident plane wave (see \cite[Theorem 5.5]{CK19} {\color{xxx}and} \cite{EH2019}). The proof carries over to the case of crack scatting as shown in the sequel.

There exist many numerical approaches for inverse scattering problems such as iterative solution method, decomposition method and sampling method (cf. \cite{CK19}).
Recently, the so-called extended sampling method has been proposed in \cite{LS2018} to determine the location and approximate the support of unknown scatterers from {\color{xxx} far-field data} generated by one incident plane wave.
As a variant of the classical linear sampling method \cite{Cakoni14,Colton1996}, extended sampling method is based on the indicator function $z\mapsto ||g_z||_{L^2(\Sp)}$, where the function $g_z$ with $z\in \R^2$ {\color{xxx}is a regularized solution to} the integral equation
\be\label{eq0}
\int_{\Sp}u_{B_R(z)}^\infty(\hat x,d)g_z(d)ds(d)=U(\hat x),\quad\hat x\in\Sp.
\en
The right hand side $U(\hat x)$ of \eqref{eq0} denotes the measurement far-field data to the unknown scatterer and $u_{B_R(z)}^\infty(\hat x,d)$ is the far-field pattern to the sound-soft disk $B_R(z):=\{x\in \R^2:|x-z|<R\}$ incited by the incident plane wave with the direction $d\in\Sp$.
It has been proved in \cite{LS2018} that the Herglotz wave function $v_{g_z^\varepsilon}$ with kernel given by the regularized solution $g_z^\varepsilon$ to (\ref{eq0}) converges in $H^1(B_R(z))$ as $\varepsilon\ra0$ if $D\subset B_R(z)$ and blows up in the norm of $H^1(B_R(z))$ as $\varepsilon\ra0$ if $D\cap B_R(z)=\emptyset$.
$B_R(z)$ can be viewed as a ``test domain'' and the idea of ``range test'' can also be found in \cite{Potthast_2003}.
Combining the idea of ``test domain'' and the classical factorization method (see \cite{Kirsch08}), a variant to factorization method with one plane wave has been proposed in \cite{EH2019,mgq}.
We will apply this method to the detection of an unknown crack based on the far-field pattern corresponding to one incident plane wave or point source.

This paper is organized as follows. Uniqueness results for inverse piecewise linear crack scattering problems with a single incident wave are established in Section \ref{sec2}.
The contrast sampling method will be proposed in Section \ref{sec-indi} to initially determine the detection area.
We introduce the one-wave factorization method in Section \ref{sec3}, which will be used to roughly recover the shape and location of the unknown crack.
To get a more accurately reconstruction, we then employ Newton's iteration method in Section \ref{iteration_method}.
Numerical examples are illustrated in Section \ref{sec4}.
%Finally, a conclusion will be given in Section \ref{sec5}.

\section{Uniqueness results}\label{sec2}
\setcounter{equation}{0}

In this section, we will prove some uniqueness results for inverse crack scattering problems with a single plane wave or point source.
Denote by $u^s_j$, $u_j$ and $u^\infty_j$ the scattered field, total field and far-field pattern, respectively, associated with the crack $\G_j$ $(j=1,2)$ and corresponding to the incident plane wave $u^i$. The corresponding notations to
the incident point source $w^i$ will be denoted by $w^s_j$, $w_j$ and $w^\infty_j$, respectively. Below we state and prove
the uniqueness results.

\begin{theorem}\label{thm2.1}
Assume $\G_j$ is a sound-soft crack such that $\G_j\subset\pa\Om_j$, where $\pa\Om_j$ denotes the boundary of some convex polygon $\Om_j\subset\R^2$, {\color{xxx}$j=1,2$,} as shown in Figure \ref{Case1}.
%{\color{xxx}Let the wave number be fixed.}

(i) Let $d_0\in\Sp$ be an arbitrary fixed {\color{xxx}incident} direction. If the far-field patterns satisfy
\be\label{eq8}
u_1^\infty(\hat x,d_0)=u_2^\infty(\hat x,d_0)\quad \mbox{for all}\quad \hat x\in\Sp_0,
\en
where $\Sp_0$ is an open subset of $\Sp$, then $\G_1=\G_2$.

(ii) Let $y_0\in G$ be an arbitrarily fixed {\color{xxx}point}, where $G$ denotes the unbounded component of the complement of $\ov{\G_1\cup\G_2}$.
If the far-field patterns satisfy
\be\label{xeq8}
w_1^\infty(\hat x,y_0)=w_2^\infty(\hat x,y_0)\quad\mbox{for all}\quad \hat x\in\Sp_0,
\en
then $\G_1=\G_2$.
Here $\Sp_0$ is again an open subset of $\Sp$.
\end{theorem}

\begin{proof}
(i).
Assume to the contrary that $\G_1\neq\G_2$.
We deduce from (\ref{eq8}) by analyticity that
\ben%\label{eq8+}
u_1^\infty(\hat x,d_0)=u_2^\infty(\hat x,d_0),\quad\hat x\in\Sp.
\enn
From Rellich's lemma \cite[Theorem 2.14]{CK19}, we deduce that $u_1^s(\cdot,d_0)=u_2^s(\cdot,d_0)$ in $G$ and thus
\be\label{eq10}
u_1(x,d_0)=u_2(x,d_0),\quad x\in G.
\en
Without loss of generality, we may assume that $\G_2\ba\ov{\G_1}$ is nonempty.
Noting that $\Om_2$ is convex, we can find a line segment $\G_0\subset\G_2\ba\ov{\G_1}$ such that $\Om_2$ is located in one of the half-space, denoted by $\R^2_+$, divided by the infinity straight line containing $\G_0$.

We claim that $u_1$ can be analytically extended as an odd function with respect to an infinite straight line.
To show this, we consider the following two cases:
\begin{figure}[htbp]
  \centering
  \subfigure[Case (a)]{
    \includegraphics[width=0.45\textwidth]{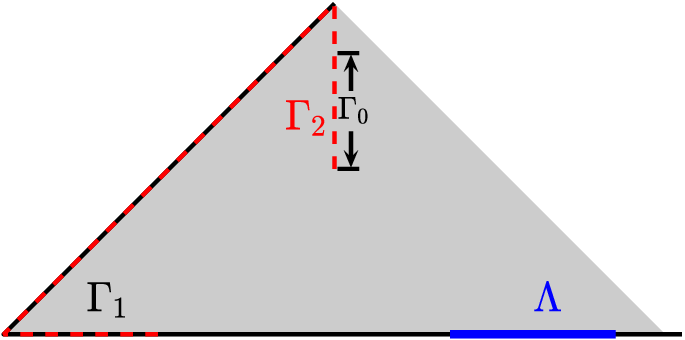}}
  \subfigure[Case (b)]{
  \includegraphics[width=0.45\textwidth]{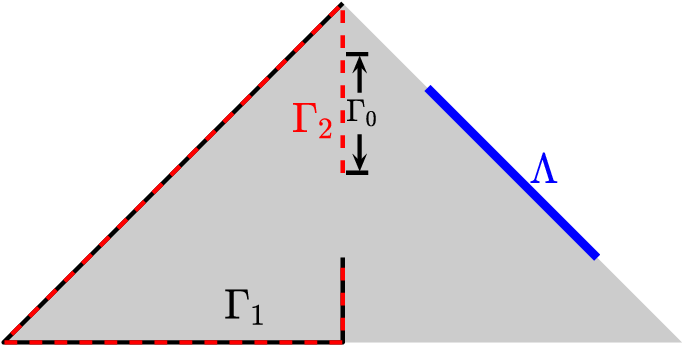}}
  \caption{An example for geometry of two cracks. Black '-': $\G_1$; red '- -': $\G_2$.}\label{Case1}
\end{figure}

\textbf{Case (a)}: $\G_1\ba\ov{\R^2_+}\neq\emptyset$ (Figure \ref{Case1} (a)).

In this case there exists a line segment $\Lambda\subset\G_1\ba\ov{\R^2_+}$. Since
$\Om_2\subset\R^2_+$, we have $\Lambda\cap\Om_2=\emptyset$.
Due to the sound-soft boundary condition on $\G_1$, (\ref{eq10}) implies $u_2(x,d_0)=u_1(x,d_0)=0$ for $x\in\Lambda\subset\G_1$.
By the analyticity of $u_2$ in $\R^2_-:=\R^2\ba\ov{\R^2_+}$, $u_1(\cdot,d_0)=u_2(\cdot,d_0)$ vanishes identically on the half line extending $\Lambda$ to infinity in {\color{xxx}$\R^2_-$}.

\textbf{Case (b)}: $\G_1\subset\R^2_+$ (Figure \ref{Case1} (b)).

Then by the sound-soft boundary condition on $\G_0\subset\G_2$, we deduce from (\ref{eq10}) that $u_1(x,d_0)=u_2(x,d_0)=0$ for $x\in\G_0$.
By the analyticity of $u_1(\cdot,d_0)$ near $\Lambda$ and the reflected principle for the Helmholtz equation (see \cite[Theorem 2.18]{Kirsch08} and \cite{ZZ13}), $u_1(\cdot,d_0)$ must be an odd function with respect to $\G_0$ in the neighborhood of $\G_0$.
Noting that $\G_1\subset\R^2_+$, we can find a line segment $\Lambda$ as a subset of the reflection of $\G_1$ with respect to $\G_0$ (see Figure \ref{Case1} (b)).
By the sound-soft boundary condition on $\G_1$ we have $u_1(x,d_0)=0$ for $x\in\G_1$.
This implies $u_1(x,d_0)=0$ for $x\in\Lambda$, since $u_1(\cdot,d_0)$ is analytic and odd with respect to $\G_0$.
Analogously to \textbf{Case (a)}, $u_1(\cdot,d_0)$ can be analytically extended along the half line that contains $\Lambda$ and lies in $\R^2_-$.

Combining the above two cases, we have proved that $u^i(\cdot,d_0)+u_1^s(\cdot,d_0)=u_1(\cdot,d_0)=0$ on the half line containing $\Lambda$ and extending to infinity in $\R^2_-$.
This is a contradiction, because $u_1^s(x,d_0)\ra0$ as $|x|\ra\infty$ by (\ref{eq13}) and $|u^i(x,d_0)|=1$ for $x\in\R^2$.

(ii).
Assume to the contrary that $\G_1\!\neq\!\G_2$.
We deduce from (\ref{xeq8}) by analyticity that
\ben%\label{818}
w_1^\infty(\hat x,y_0)=w_2^\infty(\hat x,y_0)
\enn
for all $\hat x\in\Sp$. By Rellich's lemma \cite[Theorem 2.14]{CK19}, we have $w_1^s(\cdot,y_0)=w_2^s(\cdot,y_0)$ in $G$ and thus
\ben%\label{818-1}
w_1(x,y_0)=w_2(x,y_0)
\enn
for all $x\in G\ba\{y_0\}$.
Analogously to the proof of assertion (i), we know that $w_1(\cdot,y_0)$ must be an odd function in a symmetric subdomain of $\R^2\backslash\overline{\Gamma_2}$ with respect to some straight line $\widetilde\Lambda$ (e.g., Figure \ref{Case1}).
Hence, by the sound-soft boundary condition on $\G_1$ we have
\be\label{818-2}
w_1(x,y_0)=0,\quad x\in\G_1\cup\widetilde\Lambda.
\en
Denote the reflected point of $y_0$ with respect to $\widetilde\Lambda$ by $y'_0$.
Due to the singularity of the fundamental solution, $w_1(x,y_0)$ is singular at $x=y_0$.
Thus, $w_1(x,y_0)$ is singular at $x=y'_0$, since $w_1$ is odd.
However, this is impossible because, $w_1(y'_0,y_0)=0$ by (\ref{818-2}) if $y'_0\in\G_1\cup\widetilde\Lambda$ and $w_1(\cdot,y_0)$ is analytic in the neighborhood of $y'_0$ if $y'_0\notin\G_1\cup\widetilde\Lambda$.
\end{proof}

\begin{remark}\label{rem2.2}
(i) Theorem \ref{thm2.1} remains valid even if  $\Om_1$ and $\Om_2$ are non-convex polygons.
By the sound-soft boundary condition, one can apply the reflected principle finitely many times to find a {\color{xxx}half line} $\Lambda$ such that one of the total fields {\color{xxx}corresponding to $\Gamma_1$ and $\Gamma_2$} vanishes on it
and satisfies the Helmholtz equation in the neighborhood.

(ii) The proof of Theorem \ref{thm2.1} implies that the total field cannot be {\color{xxx}analytically extended across crack tips or interior corners of a} piecewise linear {\color{xxx}sound-soft} crack.
\end{remark}

\section{Contrast sampling method to determine a rough location}\label{sec-indi}
\setcounter{equation}{0}

In this section, we introduce a contrast sampling method to get a rough location of the unknown crack.
This method is based on a comparison between the far-field data {\color{xxx}of the target crack and the far-field patterns of test point sources/test scatterers}.
To show the method, we assume for a while that the shape of the target crack $\G$ is known {\color{xxx}a priori.
%Let the target crack $\G$ satisfy the assumptions of Theorem \ref{thm2.1}.
For $a\in\R^2$} denote the shifted crack by {\color{xxx}$\G_a:=\{x+a\in\R^2:x\in\G\}$}. Define the indicator function
\be\label{comp0}
I_{crack}({\color{xxx}a}): =\left\{\int_{\Sp_0}|U_{\G}^\infty(\hat x)-U_{{\color{xxx}\G_a}}^\infty(\hat x)|^2ds(\hat x)\right\}^{-1},
\en
where $U_{\G}^\infty(\hat x)$ and $U_{{\color{xxx}\G_a}}^\infty(\hat x)$ are the far-field patterns corresponding to an incident plane wave $u^i(x,d_0)$ for an arbitrarily fixed $d_0\in\Sp$ or an incident point source $w^i(x,y_0)$ for an arbitrarily fixed $y_0\in\R^2\ba\ov{(\G\cup{\color{xxx}\G_a})}$, associated with the crack $\G$ and ${\color{xxx}\G_a}$, respectively.
Here, again {\color{xxx}$\Sp_0$} is an open subset of the unit circle $\Sp$.
By Theorem \ref{thm2.1}, the indicator {\color{xxx}$I_{crack}(a)$} is well defined and positive for all {\color{xxx}$a\in\R^2\ba\{0\}$} and $I_{crack}({\color{xxx}a})\ra\infty$ as {\color{xxx}$|a|\ra0$}.
Noting that ${\color{xxx}\G_a}=\G$ provided {\color{xxx}$|a|=0$}, we conclude that the location of the target crack can be recovered by plotting the indicator function {\color{xxx}$a\mapsto I_{crack}(a)$}.

If the shape of the target crack is unknown,
we replace $U_{{\color{xxx}\G_a}}^\infty(\hat x)$ in (\ref{comp0}) by the far-field pattern of a test disk incited by the same incident wave or the far-field pattern of a test point source.
This leads to our contrast sampling method.

\subsection{Comparison with point sources}

Replacing $U_{{\color{xxx}\G_a}}^\infty(\hat x)$ in (\ref{comp0}) by the far-field pattern of a test point source, we obtain the indicator function
\be\label{comparison2}
I_{ps}(P;\tau):=\left\{\int_{\Sp_0}\left|U_\G^\infty(\hat x)-\tau{\color{xxx}\Phi_k^\infty(\hat x,P)}\right|^2ds(\hat x)\right\}^{-1},
\en
where {\color{xxx}$\Phi_k^\infty(\hat x,P)$ is the far-field pattern of the point source $\Phi_k(x,P)$ located at $P\in\R^2$}, and {\color{xxx}$\tau\in\C\ba\{0\}$} is the scattering strength.
Here, the wave number {\color{xxx}$k>0$ of the point source $\Phi_k(x,P)$ is the same as the incident wave of the target crack}.
As explained in the beginning of this section, it is expected that $I_{ps}(P;\tau)$ will take a large value as $P$ is getting closer to $\G$ and take relatively small values when $P$ moves away from $\G$.

\subsection{Comparison with disks}

Suppose $B_R(P)$ is a disk centered at $P\in\R^2$ with radius $R>0$.
Replacing $U_{{\color{xxx}\G_a}}^\infty(\hat x)$ in (\ref{comp0}) by the far-field pattern $U_{B_R(P)}^\infty(\hat x)$ of the sound-soft disk $B_R(P)$ incited by the same incident wave, we obtain the indicator function
\be\label{comparison1}
I_{disk}(P;R):=\left\{\int_{\Sp_0}\left|U_\G^\infty(\hat x)-U_{B_R(P)}^\infty(\hat x)\right|^2ds(\hat x)\right\}^{-1}.
\en
As explained in the beginning of this section, it is expected that $I_{disk}(P;R)$ will take a large value as $P$ is close to $\G$ and take relatively small values as $P$ moves away from $\G$.
We guess the indicator function $I_{disk}(P;R)$ defined by (\ref{comparison1}) remains valid if the boundary condition of $B_R(P)$ is replaced by other boundary conditions or even the disk is replaced by other scatterers, like an inhomogeneous medium.

In the numerical implementation of contrast sampling method given by (\ref{comparison1}), we need to solve a forward scattering problem for different $P$ or $R$.
For incident plane waves, the far-field pattern corresponding to a disk takes the explicit series form
\be\label{eig-dir}
  u_{B_R(P)}^\infty(\hat x,d_0)=-\sqrt{\frac{2}{k\pi}}e^{-i\frac{\pi}4}\sum_{n\in\Z}\frac{ J_n(kR)}{H_n^{(1)}(kR)}e^{in(\theta_{ x}-\theta_d)}e^{ikP\cdot(d-\hat x)},\quad\hat x\in\Sp,
\en
if $B_R(P)$ is sound-soft, and
\begin{equation}\label{eig-imp}
  u_{B_R(P)}^\infty(\hat x,d_0)=-\sqrt{\frac{2}{k\pi}}e^{-i\frac{\pi}4}\sum_{n\in\Z}\frac{kJ'_n(kR)+\eta J_n(kR)}{kH_n^{(1)\prime}(kR)+\eta H_n^{(1)}(kR)}e^{in(\theta_{ x}-\theta_d)}e^{ikP\cdot(d-\hat x)},\quad\hat x\in\Sp,
\end{equation}
if the impedance boundary condition $\pa_\nu u+\eta u=0$ is {\color{xxx}imposed} on $\pa B_R(P)$ with the constant impedance coefficient $\eta\in \C$.
Here, $\hat x =(\cos \theta_x, \sin \theta_x)$, $d_0=(\cos \theta_d, \sin \theta_d)$, $J_n$ is the Bessel function of order $n$, and $H^{(1)}_n$ is the Hankel function of the first kind of order $n$. The expressions
(\ref{eig-dir}) and (\ref{eig-imp}) can be deduced from the expansion of the scattered field in terms of spherical wave functions together with the translation property (see e.g., \cite[(2.49) and (5.3)]{CK19} and \cite{mgq}).

\section{The one-wave factorization method}\label{sec3}
\setcounter{equation}{0}

Using the contrast sampling method in the previous section, the rough location of target crack can be determined by far-field data of a single wave.
In this section, we will employ the one-wave factorization method introduced in \cite{mgq} for getting a precise information on the location and shape of the target crack.
We first review {\color{xxx}the classical factorization method (cf. \cite{Kirsch08})}.

\subsection{Preliminary results from the factorization method}\label{sub4.1}

Suppose $D\subset\R^2$ is a bounded obstacle with its boundary $\pa D \in C^2$ such that $\R^2\backslash\overline{D}$ is connected.
Define the data-to-pattern operator corresponding to $D$ by $G_Df=u_D^\infty$ where $u_D^\infty$ is the far-field pattern to the solution $u_D^s$ of the following boundary value problem
\be\label{eq17}
\Delta u^s_D+k^2u^s_D=0 & \text{in}\;\R^2\ba\ov{D},\\ \label{eq18}
\mathscr Bu^s_D=f & \text{on}\;\pa D,\\ \label{eq19}
\lim_{r\ra\infty}\sqrt{r}\left(\frac{\pa u^s_D}{\pa r}-iku^s_D\right)=0, & r=|x|,
\en
where the boundary condition $\mathscr B$ on $\pa D$ depends on the physical property of $D$:
\ben%\label{eq21}
\mathscr Bu=\begin{cases}
u, & \text{if}\;D\;\text{is a sound-soft obstacle},\\
\frac{\pa u}{\pa\nu}+\eta u, & \text{if}\;D\;\text{is an impedance obstacle}.
\end{cases}
\enn
Here, $\nu$ denotes the unit outward normal to $\pa D$ and $\eta\!\in\!C(\pa D)$ is the impedance coefficient satisfying ${\rm Im}\,\eta\geq0$ on $\pa D$.
In particular, $D$ is called a sound-hard obstacle if $\eta=0$ on $\pa D$.
The existence of a unique solution to above exterior boundary value problem can be established either by integral equation method \cite[Chapter 3]{CK83} or by variational method \cite[Section 5.3]{Cakoni14}.
Therefore, the data-to-pattern operator $G_Df\in L^2(\Sp)$ is well-defined for all $f\!\in\!H^{1/2}(\pa D)$ if $D$ is a sound-soft obstacle and for all $f\!\in\!H^{-1/2}(\pa D)$ if $D$ is an impedance obstacle.

In particular, denote by $u_D^\infty(\hat x,d)$ the far-field pattern of the solution $u^s_D$ to (\ref{eq17})--(\ref{eq19}) with $f=-\mathscr Bu^i(\cdot,d)$ with the incident plane wave $u^i$ given by (\ref{eq-inc}).
{\color{xxx}Define the far-field operator $F_D:L^2(\s) \to L^2(\s)$ by}
\be\label{eq26}
(F_D g) (\hat x) := \int_{\s} u_D^{\infty}(\hat x,d) g(d)\,ds(d),\quad\hat x \in \Sp.
\en
The {\color{xxx}classical} factorization method is mainly based on the following {\color{xxx}theorem, which follows easily from \cite[Theorems 1.21 and 2.15]{Kirsch08} and \cite[Lemmas 3.3 and 3.4]{XZZ}.}
\begin{theorem}%\label{thm-ri}
If $k^2$ is not an eigenvalue of $-\Delta$ in $D$ with corresponding boundary condition, i.e., the following interior boundary value problem
\ben
\Delta v+k^2v=0 && \text{in }D,\\
\mathscr Bv=0 && \text{on }\pa D,
\enn
has only the trivial solution $v=0$ in $D$.
Then we have
\be\label{ri}
{\color{xxx}{\rm Ran}\,F_{D,\#}^{1/2}={\rm Ran}\,G_D,}
\en
where {\color{xxx}${\rm Ran}\,F_{D,\#}^{1/2}$ and ${\rm Ran}\,G_D$ denote the ranges of the operators $F_{D,\#}^{1/2}$ and $G_D$, respectively, and} $F_{D,\#}:=|{\rm Re}\,F_D|+|{\rm Im}\,F_D|$ with ${\rm Re}\,F_D:=(F_D+F_D^*)/2$ and {\color{xxx}${\rm Im}\,F_D:=(F_D-F_D^*)/(2i)$.}
\end{theorem}
%\begin{remark}%\label{rem240323}
%If $D$ is an impedance obstacle with ${\rm Im}\,\eta>0$ in an open subset of $\pa D$, then $k^2>0$ is not an impedance eigenvalue of $-\Delta$ in $D$ (see \cite[Theorem 8.2]{Cakoni14}).
%\end{remark}

\subsection{The factorization method with a single wave}
{\color{xxx}Let $D$ be given as in Subsection \ref{sub4.1}. In this subsection $D$ will play the role of testing scatterers for detecting the location and shape of a piecewise linear crack, as shown by the following theorem on the one-wave factorization method.}

\begin{theorem}\label{t1}
Assume $\G$ is a sound-soft crack such that there exists a polygon $\Om$ whose boundary $\pa\Om$ satisfies $\G\subset\pa\Om$.
Suppose $U_\G^\infty$ is the far-field pattern to $\G$ corresponding to the plane wave (\ref{eq-inc}) with an arbitrarily fixed $d_0\in\Sp$ or the point source (\ref{ps}) with an arbitrarily fixed $y_0\in\R^2\ba\ov{\G}$.

(a) If $\G\subset D$ (see Figure \ref{insideD_test} (a)), then $U_\G^\infty\in{\rm Ran}\,G_D$.

{\color{xxx}(b) If $\Gamma\not\subset D$ and $D$ is convex (see Figure \ref{insideD_test} (b)), then $U_\G^\infty\notin{\rm Ran}\,G_D$.

(c) If $\G\cap D=\emptyset$ (see Figure \ref{insideD_test} (c)), then $U_\G^\infty\notin{\rm Ran}\,G_D$.}
\end{theorem}

\begin{proof}
Let $U_\G^s$ and $U_\G$ denote the scattered field and total field of the crack $\G$ corresponding to the far-field pattern $U_\G^\infty$, respectively.

(a). Since $\G\subset D$, we can set $f=\mathscr B U_\G^s$ in (\ref{eq18}).
The uniqueness of the exterior boundary value problem (\ref{eq17})--(\ref{eq19}) implies $u^s_D=U_\G^s$ in $\R^2\ba\ov{D}$.
Therefore, $u_D^\infty(\hat x)=U_\G^\infty(\hat x)$ for all $\hat x\in\Sp$ and thus $U_\G^\infty=G_Df$.

(b). Assume to the contrary that $U_\G^\infty\in{\rm Ran}\,G_D$.
Then there exists a boundary value $f$ on $\pa D$ such that the far-field pattern $u_D^\infty$ of the solution $u_D^s$ to the exterior boundary value problem (\ref{eq17})--(\ref{eq19}) coincides with $U_\G^\infty$.
Rellich's lemma \cite[Theorem 2.14]{CK19} implies $u^s_D=U_\G^s$ in the unbounded component of the complement of $\G\cup D$.
Noting that $\G\not\subset D$ and $D$ is convex, one can find a line segment $\Lambda_0\subset\G\ba\ov{D}$ such that the total field $U_\G$ is analytic in the neighborhood of $\Lambda_0$ and vanishes on $\Lambda_0$.
{\color{xxx}In particular, by the convexity of $D$, one can always assume that the line segment $\Lambda_0$ contains a crack tip or interior corner of $\Gamma$}.
In view of Remark \ref{rem2.2}\,(ii), proceeding as in the proof of Theorem \ref{thm2.1}, we can obtain a contradiction between the analyticity of $u_D^s$ and the singularity of $U_\Gamma^s$ at {\color{xxx}the crack tip or interior corner}.

\begin{figure}[htb]
  \centering
  \subfigure[$\Gamma\subset D$]{
  \includegraphics[width=0.2\textwidth]{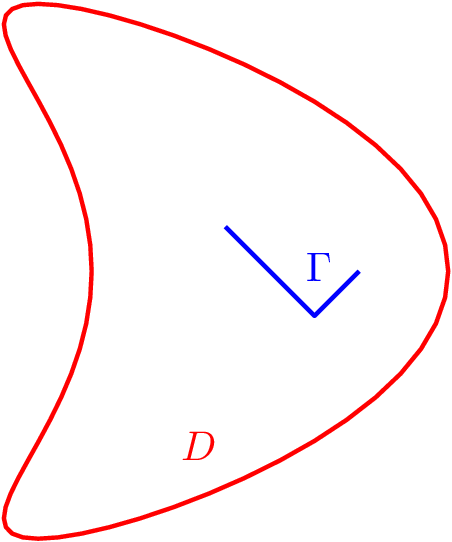}}
  \subfigure[$\Gamma\not\subset D$]{
  \includegraphics[width=0.25\textwidth]{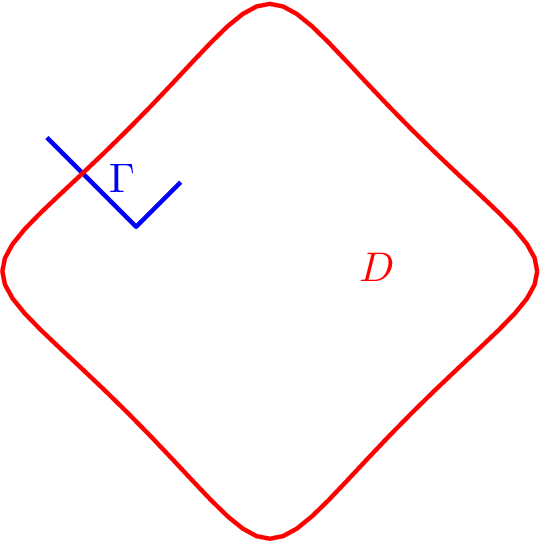}}
  \subfigure[$\Gamma\cap D=\emptyset$]{
  \includegraphics[width=0.27\textwidth]{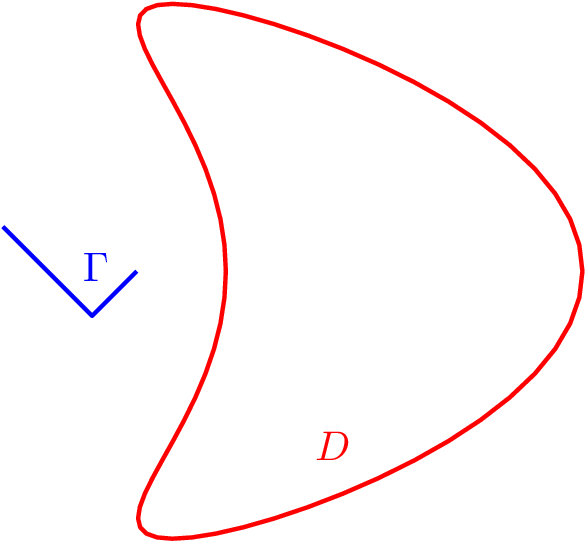}}
  \subfigure[$\G$, concave $D$]{
  \includegraphics[width=0.2\textwidth]{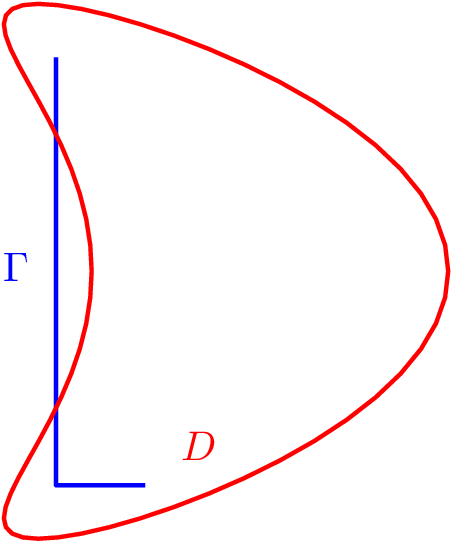}}
  %\subfigure[Counterexample]{
  %\includegraphics[width=0.18\textwidth]{fig/reflection2.eps}}
  \caption{Geometry of Theorem \ref{t1}.}\label{insideD_test}
\end{figure}

{\color{xxx}(c). Assume to the contrary that $U_\G^\infty\in{\rm Ran}\,G_D$.
Proceeding as in the proof of assertion (b), we can conclude that $U_\G^s$ can be extended as an entire solution to the Helmholtz equation.
Therefore, the radiating solution $U_\G^s$ vanishes identically.
However, this contradicts the sound-soft boundary condition on $\G$.}
\end{proof}
\begin{remark} For non-convex testing scatterers $D$, {\color{xxx}the assertions (a) and (c) of Theorem \ref{t1} remain valid}, but the assertion (b) is no longer true. In fact, if {\color{xxx}the crack tips and corner points of $\Gamma$} are all contained in {\color{xxx}a concave domain $D$ but $\Gamma\not\subset D$} (see Figure \ref{insideD_test} (d)), one can prove via analytical continuation arguments that {\color{xxx}$U_\G^\infty\in{\rm Ran}\,G_D$}.
\end{remark}

Theorem \ref{t1} implies that the inclusion relation between the crack $\G$ and the test domain $D$ can be characterized by whether the far-field pattern $U_\G^\infty$ to the crack belongs to ${\rm Ran}\,G_D$ or not.
Motivated by this, a numerical method to reconstruct the location and rough profile of target crack $\G$ from far-field data $U_\G^\infty$ of a single wave can be designed.
There are several methods to compute the range of $G_D$.
To solve the inverse crack scattering problem in a data-to-data manner, we shall get ${\rm Ran}\,G_D$ indirectly from the range identity (\ref{ri}) that requires the far-field data for $D$, instead of by directly solving the exterior boundary value problem (\ref{eq17})--(\ref{eq19}) (see \cite{Potthast_2003}).
This leads to our factorization method with a single wave:

\begin{theorem}\label{thm3.9}
Suppose that the test domain $D$ introduced in Subsection \ref{sub4.1} is convex.
Assume that $k^2>0$ is not an eigenvalue of $-\Delta$ over $D$ with the corresponding boundary condition.
Denote by $(\lambda_n,f_n)$ an eigensystem of $F_{D,\#}$ and
let $U_\G^\infty$ be the same as in Theorem \ref{t1}.
Define the indicator function by
\be\label{I-1}
I_1(D):=\left\{{\color{xxx}\sum_{n}}\frac1{\lambda_n}\left|\left(U_\G^\infty,f_n\right)_{L^2(\Sp)}\right|^2\right\}^{-1}.
\en
Then
we have $I_1(D)=0$ if $\G\not\subset D$, and $I_1(D)>0$ if $\G\subset D$.
\end{theorem}

We omit the proof Theorem \ref{thm3.9}, since it follows directly from (\ref{ri}), Theorem \ref{t1}, {\color{xxx}the fact that ${\rm Ran}\,G_D$ is dense in $L^2(\Sp)$ (see \cite[Theorem 4.8]{Cakoni14} and \cite[Theorem 3.36]{CK19})}, and Picard's theorem \cite[Theorem 4.8]{CK19}.

\begin{remark}
If $D=B_R(P)$ is a disk centered at $P$ with radius $R>0$, the above theorem was shown in \cite[Theorem 3.7]{mgq}. For circular test domains, the eigensystem of $F_{B_R(P),\#}$ can be deduced directly from \eqref{eig-dir} and \eqref{eig-imp} as follows.
\begin{itemize}
  \item If $B_R(P)$ is a sound-soft disk, then the eigensystem of $F_{B_R(P),\#}$ is given by
  \ben
  \lambda_n=\sqrt{\frac{8\pi}{k}}\left|\frac{J_n(kR)}{H_n^{(1)}(kR)}\right|,\quad f_n(\theta)=\frac1{\sqrt{2\pi}}e^{in\theta}e^{-ikP\cdot(\cos\theta,\sin\theta)}{\color{xxx},\quad n\in\Z};
  \enn
  \item If $B_R(P)$ is an impedance disk as in \eqref{eig-imp}, then the eigensystem of $F_{B_R(P),\#}$ is given by
  \ben
  \lambda_n=\sqrt{\frac{8\pi}{k}}\left|\frac{kJ'_n(kR)+\eta J_n(kR)}{kH_n^{(1)\prime}(kR)+\eta H_n^{(1)}(kR)}\right|,\quad f_n(\theta)=\frac1{\sqrt{2\pi}}e^{in\theta}e^{-ikP\cdot(\cos\theta,\sin\theta)}{\color{xxx},\quad n\in\Z}.
  \enn
\end{itemize}
\end{remark}

As is shown in Theorem \ref{thm2.1}, the uniqueness for inverse problem remains valid with limited aperture far-field data of a single wave.\,%, i.e., a far-field pattern on an open subset $\Sp_0$ of $\Sp$.
We will extend our factorization method with a single wave from full aperture case to limited aperture case below.
Following the idea of factorization method with limited aperture data (see \cite[Section 2.3]{Kirsch08}), we introduce the operator $F_{D,la}:L^2(\Sp_0)\ra L^2(\Sp_0)$:
\be\label{Fla}
(F_{D,la}g)(\hat x):=\int_{\Sp_0}u_D^\infty(\hat x,d)g(d)ds(d),\quad\hat x\in\Sp_0.
\en
By \cite[(2.49)]{Kirsch08} we have the following relation between (\ref{eq26}) and (\ref{Fla}):
\be\label{eq-3}
F_{D,la}=P_{\Sp_0}F_DP_{\Sp_0}^*,
\en
where $P_{\Sp_0}:L^2(\Sp)\ra L^2(\Sp_0)$ is the restriction operator $P_{\Sp_0}g=g|_{\Sp_0}$.
The adjoint $P_{\Sp_0}^*:L^2(\Sp_0)\ra L^2(\Sp)$ is the zero extension, i.e., $(P_{\Sp_0}^*g)(d)=g(d)$ for $d\in\Sp_0$ and $(P_{\Sp_0}^*g)(d)=0$ otherwise.
{\color{xxx}Analogously to \eqref{ri}, we conclude from \eqref{eq-3} that ${\rm Ran}\,F_{D,la,\#}^{1/2}={\rm Ran}(P_{\Sp_0}G_D)$ (see \cite[Theorems 2.9 and 2.15]{Kirsch08}). %, where $F_{D,la,\#}=|{\rm Re}\,F_{D,la}|+|{\rm Im}\,F_{D,la}|$ with ${\rm Re}\,F_{D,la}=(F_{D,la}+F_{D,la}^*)/2$ and ${\rm Im}\,F_{D,la}=(F_{D,la}-F_{D,la}^*)/(2i)$.
By analyticity, $U_\G^\infty|_{\Sp_0}\in{\rm Ran}\,(P_{\Sp_0}G_D)$ is equivalent to $U_\G^\infty\in{\rm Ran}\,G_D$.
In view of Theorem \ref{t1}, we immediately obtain the following theorem.}

\begin{theorem}\label{thm3.9-la}
{\color{xxx}Denote by $(\lambda_{la,n},f_{la,n})$ an eigensystem of $F_{D,la,\#}$.}
Under the assumptions of Theorem \ref{thm3.9}, we have $I_2(D)=0$ if $\G\not\subset D$, and $I_2(D)>0$ if $\G\subset D$.
Here, the indicator function is defined by
\be\label{211005-3}
I_2(D):=\left\{{\color{xxx}\sum_{n}}\frac1{\lambda_{la,n}}\left|\left(U_\G^\infty,f_{la,n}\right)_{L^2(\Sp_0)}\right|^2\right\}^{-1}.
\en
\end{theorem}

\begin{remark}%\label{rem+1}
Note that the far-field operator $F_D$ is compact from $L^2(\Sp^2)$ to itself, since $F_D$ is an integral operator with the smooth integral kernel $u_D^\infty(\hat x,d)$ (see \cite[Theorem 1.7]{Kirsch08}).
It is easily seen that $F_D^*$, ${\rm Re}\,F_D$, ${\rm Im}\,F_D$, $|{\rm Re}\,F_D|$, $|{\rm Im}\,F_D|$, $F_{D,\#}$ and $F_{D,la,\#}$, are all compact from $L^2(\Sp^2)$ to itself.
Therefore, $\lambda_n,\lambda_{la,n}\ra0$ as $n\ra\infty$ and it is not stable to calculate (\ref{I-1}) and (\ref{211005-3}), especially when the far-field pattern $U_\G^\infty(\hat x)$ is noise-polluted by \be\label{noi}
U_{\G,\delta}^\infty(\hat x)=U_\G^\infty(\hat x)+\delta(\zeta_1+i\zeta_2)|U_\G^\infty(\hat x)|,
\en
with $\delta$ denoting the noise ratio and $\zeta_1,\zeta_2$ being the {\color{xxx}uniformly} distributed random numbers in $[-1,1]$.
For a more stable numerical result, we apply the Tikhonov regularization (see \cite[Section 4.4]{CK19}) to (\ref{I-1}) and (\ref{211005-3}) to obtain
\be\label{al3.1}
I_1(D)\!\!\!\!&\approx&\!\!\!\!{\color{xxx}\widetilde I_1(D):=}\left\{{\color{xxx}\sum_{n}}\frac{\lambda_n}{(\alpha+\lambda_n)^2}\left|\left(U_{\G,\delta}^\infty,f_n\right)_{L^2(\Sp)}\right|^2\right\}^{-1},\\
\label{211012-2}
I_2(D)\!\!\!\!&\approx&\!\!\!\!{\color{xxx}\widetilde I_2(D):=}\left\{{\color{xxx}\sum_{n}}\frac{\lambda_{la,n}}{(\alpha+\lambda_{la,n})^2}\left|\left(U_{\G,\delta}^\infty,f_{la,n}\right)_{L^2(\Sp_0)}\right|^2\right\}^{-1},
\en
respectively, where {\color{xxx}the regularization parameter $\alpha>0$ is appropriately chosen.}
\end{remark}

With this method, the location and profile of the target crack can be recovered by selecting appropriate sampling domains, such as disks with different centers and radii near the rough location of the crack given by the contrast sampling method introduced in Section \ref{sec-indi}.
Theoretically, the location and convex hull of an unknown scatterer can be recovered by the factorization method with a single wave (see \cite[Remark 3.10]{mgq}).
It should be pointed out that the above method can be extended to the case when the obstacle $D$ is replaced by {\color{xxx}an inhomogeneous medium whose contrast function has a compact support $D$, provided $k^2$ is not a corresponding interior transmission eigenvalue.}

\section{A more precise result by Newton's iteration method}\label{iteration_method}
\setcounter{equation}{0}

For a more precise numerical result, we can apply Newton's iteration method whose initial guess is given by the method introduced in the previous section.
On the other hand, a proper initial guess also improves the behavior and result of iteration method.
{\color{xxx}For details on the iteration method, we refer the reader to \cite{Bochniak2002,Kress1995}.

We begin with the numerical simulation of forward crack scattering. A piecewise} linear crack with two tips and interior corners given in order by $\{P_{\ell}:=(P_{\ell,1},P_{\ell,2})\}_{\ell=0,\cdots,N}$ possesses a parametric representation of the form $x(t):=(x_1(t),x_2(t))$, $0\leq t\leq2\pi$, where, for $j=1,2$,
\be\label{230912-1}
x_j(t)\!=\!\left(\frac{2(\ell\!+\!1)\pi}{N}\!-\!t\right)P_{\ell,j}\!+\!\left(t\!-\!\frac{2\ell\pi}{N}\right)P_{\ell+1,j},\;t\in\left[\frac{2\ell\pi}{N},\frac{2(\ell\!+\!1)\pi}{N}\right),\ell\!=\!0,\cdots,N\!-\!1.
\en
For a more precise numerical simulation of the far-field pattern and scattered field, we make use of a graded mesh rather than a uniform mesh (see \cite[Section 3.6]{CK19}).
To introduce the graded mesh, we choose $n\in\N$ such that $n/N\in\Z$.
For simplicity let there be $n/N$ knots on each smooth segment% (the number of knots on each smooth segment can, in general, be arbitrary)
.
The boundary knots of the graded mesh are given by $x(t_j)$ with
\ben
t_j=w(s_j),\quad s_j=\frac{\pi}{2n}+\frac{j\pi}{n},\quad j=0,1,\cdots,2n-1,
\enn
where
\ben
&&w(s)={\color{xxx}\tilde w(Ns-2\ell\pi)},\quad s\in\left[\frac{2\ell\pi}{N},\frac{2(\ell+1)\pi}{N}\right),\ell=0,\cdots,N-1,\\
&&\tilde w(s)=2\pi\frac{[v(s)]^p}{[v(s)]^p+[v(2\pi-s)]^p},\quad 0\leq s\leq2\pi,\\
&&v(s)=\left(\frac1p-\frac12\right)\left(\frac{\pi-s}\pi\right)^3+\frac{s-\pi}{p\pi}+\frac12,\quad p\geq2.
\enn
According to \cite[Section 8.7]{Cakoni14}, the scattered field and its far-field pattern to this crack can be represented as
\be\label{211013-1}
U_\G^s(x)\!\!\!\!&=&\!\!\!\!\int_{\G}\Phi_k(x,y)\varphi(y)ds(y),\quad x\in\R^2\ba\ov{\G},\\ \label{211013-1far}
U_\G^\infty(\hat x)\!\!\!\!&=&\!\!\!\!\frac{e^{i\frac\pi4}}{\sqrt{8k\pi}}\int_{\G}e^{-ik\hat x\cdot y}\varphi(y)ds(y),\quad\hat x\in\Sp,
\en
where the density $\varphi\in\widetilde H^{-1/2}(\G):=\{u\in H^{-1/2}(\pa\Om):{\rm supp}\,u\subset\ov{\G}\}$ solves the boundary integral equation
\ben
2\int_{\G}\Phi_k(x,y)\varphi(y)ds(y)=-2U_\Gamma^i(x),\quad x\in\G.
\enn
Following \cite[Section 3.6]{CK19}, the above boundary integral equation can be approximated by the following linear system
\be\label{0916-1}
{\color{xxx}M_{corner}}W\Psi=F,
\en
where
\ben
{\color{xxx}M_{corner}}\!\!\!\!&:=&\!\!\!\!\left(R_j(s_i)\widetilde M_1(s_i,s_j)+\frac{\pi}n\widetilde M_2(s_i,s_j)\right)_{i,j=0,1,\cdots,2n-1},\\
R_j(s)\!\!\!\!&:=&\!\!\!\!-\frac{2\pi}n\sum_{m=1}^{2n-1}\frac1m\cos m(s-s_j)-\frac{\pi}{n^2}\cos n(s-s_j),\;j=0,1,\cdots,2n-1,\\
\widetilde M_1(s,\sigma)\!\!\!\!&:=&\!\!\!\!M_1(s,\sigma),\\
\widetilde M_2(s,\sigma)\!\!\!\!&:=&\!\!\!\!\begin{cases}
M_2(s,\sigma), & s\neq\sigma,\\
M_2(w(s),w(s))+2M_1(w(s),w(s))\ln w'(s),  & s=\sigma,
\end{cases}\\
M_1(t,\tau)\!\!\!\!&:=&\!\!\!\!-\frac1{2\pi}J_0(k|x(t)-x(\tau)|)|x'(\tau)|,\\
M_2(t,\tau)\!\!\!\!&:=&\!\!\!\!\begin{cases}
M(t,\tau)-M_1(t,\tau)\ln\left(4\sin^2\frac{t-\tau}2\right), & t\neq\tau,\\
\left\{\frac i2-\frac C\pi-\frac1\pi\ln\left(\frac k2|x'(t)|\right)\right\}|x'(t)|, & t=\tau,
\end{cases}\\
M(t,\tau)\!\!\!\!&:=&\!\!\!\!\frac{i}2H_0^{(1)}(k|x(t)-x(\tau)|)|x'(\tau)|,\\
W\!\!\!\!&:=&\!\!\!\!{\rm diag}((w'(s_0),w'(s_1),\cdots,w'(s_{2n-1}))),\\
\Psi\!\!\!\!&:=&\!\!\!\!(\varphi(x(t_0)),\varphi(x(t_1)),\cdots,\varphi(x(t_{2n-1})))^\top,\\
F\!\!\!\!&:=&\!\!\!\!-2(U_\G^i(x(t_0)),U_\G^i(x(t_1)),\cdots,U_\G^i(x(t_{2n-1})))^\top,
\enn
where $C$ denotes the Euler's constant. We have $U_\G^i(x)=u^i(x,d_0)$ if the incident field is the plane wave (\ref{eq-inc}) with an arbitrarily fixed $d_0\in\Sp$ and $U_\G^i(x)=w^i(x,y_0)$ if the incident field is the point source (\ref{ps}) with an arbitrarily fixed $y_0\in\R^2\ba\ov{\G}$.
With the above notations in discrete form, the scattered field \eqref{211013-1} and its far-field pattern \eqref{211013-1far} can be approximated by
\ben%\label{211013-3}
U_\G^s(x)\!\!\!\!&\approx&\!\!\!\!\frac\pi n\left(\begin{array}{cccc}
\Phi_k(x,x(t_0)) & \Phi_k(x,x(t_1))& \cdots & \Phi_k(x,x(t_{2n-1}))
\end{array}\right)D_{\Gamma}W\Psi,\quad x\in\R^2\ba\ov{\G},\\ %\label{211013-4}
U_\G^\infty(\hat x)\!\!\!\!&\approx&\!\!\!\!\frac\pi n\frac{e^{i\frac\pi4}}{\sqrt{8k\pi}}\left(\begin{array}{cccc}
e^{-ik\hat x\cdot x(t_0)} & e^{-ik\hat x\cdot x(t_1)}& \cdots & e^{-ik\hat x\cdot x(t_{2n-1})}
\end{array}\right)D_{\Gamma}W\Psi,\quad\hat x\in\Sp,
\enn
where 
\be\label{dgamma}
D_{\Gamma}={\rm diag}(|x'(t_0)|,|x'(t_1)|,\cdots,|x'(t_{2n-1})|).
\en
\begin{remark}\label{0916rem}
Since $w'(s_j)$ takes a very small value if the knot $x(w(s_j))$ is close to the tips or interior corners of the crack, it is not stable to calculate $\Psi$ from (\ref{0916-1}).
From the above approximation for the scattered field and its far-field pattern, we know $W\Psi$ can be viewed as an unknown vector and it is sufficient to calculate $W\Psi$ from (\ref{0916-1}).
\end{remark}

To introduce the iteration method, we consider a crack $\G$ represented by $h(t)\!\in\!\R^2$ for $t\!\in\!(0,2\pi)$.
Define the far-field mapping $\mathcal F$ by
\be\label{0916-2}
(\mathcal Fh)(\hat x)=U_\G^\infty(\hat x),\quad\hat x\in\Sp_0,
\en
where $\Sp_0$ is an open subset of $\Sp$.
The inverse problem with the limited aperture far-field data can be formulated as the operator equation (\ref{0916-2}) for finding $h$, which is going to be approximately solved by Newton's iteration method.
Precisely, for a proper initial guess $h_0$ we compute
\ben
h_{n+1}=h_n+q_n,\quad n=0,1,\cdots,
\enn
where $q_n$ solves the linearized equation of (\ref{0916-2}):
\be\label{0916-7}
\mathcal F h_n+\mathcal F'_{h_n}q_n=U_\G^\infty\quad\text{on}\;\Sp_0.
\en
The Fr\'echet derivative in (\ref{0916-7}) is defined by
\ben
\mathcal F'_{h}q=\lim_{t\ra0}\frac{\mathcal F(h+tq)-\mathcal F h}{t}.
\enn
According to \cite[(28)]{Bochniak2002} and \cite[Theorem 6.3]{Kress1995}, the Fr\'echet derivative is given by
\be\label{0915}
F'_hq=v^\infty|_{\Sp_0},
\en
where $v^\infty$ is the far-field pattern of the unique radiating solution $v$ to the boundary value problem
\be\label{0916-4}
\Delta v+k^2v=0 && \text{in}\;\R^2\ba\ov{\G},\\ \label{0916-5}
v_\pm=-(\nu\cdot q)\pa_\nu U_{\G,\pm} && \text{on}\;\G,
\en
where $U_\G$ denotes the total field corresponding to the crack $\G$ and the subscript $\pm$ in (\ref{0916-5}) is understood in the following sense
\ben
v_\pm(x):=\lim_{t\ra+0}v(x\pm t\nu(x)),\;\pa_\nu U_{\G,\pm}(x)=\lim_{t\ra+0}\pa_\nu U_\G(x\pm t\nu(x)),\quad x\in\G.
\enn
Due to the regularity of elliptic equations, $F'_h$ is compact and (\ref{0916-7}) is ill-posed.
For a more stable numerical implementation, we may apply the Tikhonov regularization scheme to obtain
\be\label{211014-1}
q_n\approx(\alpha I+[\mathcal F'_{h_n}]^*\mathcal F'_{h_n})^{-1}[\mathcal F'_{h_n}]^*(U_\G^\infty|_{\Sp_0}-\mathcal F h_n),
\en
where the regularization parameter $\alpha\!>\!0$ is appropriately chosen.

There are two difficulties in the numerical implementation of Newton's iteration method:
\begin{itemize}
\item {\bf Problem 1}. It is difficult to calculate the values of $\pa_\nu U_{\G,\pm}$ near the tips and interior corners of the crack, due to the singularities of elliptic boundary value problems in nonsmooth domains (see \cite{Grisvard});

\item {\bf Problem 2}. The operator $F'_h$ is in general not uniquely defined by \eqref{0915}, {\color{xxx}since $F'_hq=F'_h\tilde q$ provided $q\neq\tilde q$ but $\nu\cdot q=\nu\cdot\tilde q$ on a straight line crack} (see \cite[Page 605]{Bochniak2002}).
\end{itemize}

In the following two subsections we will show how we deal with the above two difficulties, respectively.

\subsection{Computation of Fr\'echet derivatives}

For {\bf Problem 1}, we can approximate the value of $\pa_\nu U_{\G,\pm}$ near the tips and interior corners of the crack in the following manner.
In view of (\ref{211013-1}) and jump relations, the normal derivatives of the total field in (\ref{0916-5}) are given by
\ben
2\pa_\nu U_{\G,\pm}(x)=2\frac{\pa}{\pa\nu(x)}\int_\G\Phi_k(x,y)\varphi(y)ds(y)\mp\varphi(x)+2\pa_\nu U_\G^i(x),\quad x\in\G,
\enn
which can be approximated by the following linear system (cf. \cite{Kress1995b})
\be\label{0916-6}
D_{\Gamma}V_\pm=(HW\mp D_{\Gamma})\Psi+D_{\Gamma}Y,
\en
where
\ben
V_\pm\!\!\!\!&=&\!\!\!\!2\left(\pa_\nu U_{\G,\pm}(x(t_0)),\pa_\nu U_{\G,\pm}(x(t_1)),\cdots,\pa_\nu U_{\G,\pm}(x(t_{2n-1}))\right)^\top,\\
H\!\!\!\!&:=&\!\!\!\!\left(R_j(s_i)H_1(w(s_i),w(s_j))+\frac{\pi}nH_2(w(s_i),w(s_j))\right)_{i,j=0,1,\cdots,2n-1},\\
H_1(t,\tau)\!\!\!\!&:=&\!\!\!\!\begin{cases}
-\frac k{2\pi}\{x'_2(t)[x_1(\tau)-x_1(t)]-x'_1(t)[x_2(\tau)-x_2(t)]\}\frac{J_1(k|x(t)-x(\tau)|)}{|x(t)-x(\tau)|}|x'(\tau)|, & t\neq\tau,\\
0, & t=\tau,
\end{cases}\\
H_2(t,\tau)\!\!\!\!&:=&\!\!\!\!\begin{cases}
H(t,\tau)-H_1(t,\tau)\ln\left(4\sin^2\frac{t-\tau}2\right), & t\neq\tau,\\
\frac1{2\pi}\frac{x'_2(t)x''_1(t)-x'_1(t)x''_2(t)}{|x'(t)|}, & t=\tau,
\end{cases}\\
H(t,\tau)\!\!\!\!&:=&\!\!\!\!\begin{cases}
\frac{ik}2\{x'_2(t)[x_1(\tau)-x_1(t)]-x'_1(t)[x_2(\tau)-x_2(t)]\}\frac{H_1^{(1)}(k|x(t)-x(\tau)|)}{|x(t)-x(\tau)|}|x'(\tau)|, & t\neq\tau,\\
\frac1{2\pi}\frac{x'_2(t)x''_1(t)-x'_1(t)x''_2(t)}{|x'(t)|}, & t=\tau,
\end{cases}\\
Y\!\!\!\!&=&\!\!\!\!2\left(\pa_\nu U_\G^i(x(t_0)),\pa_\nu U_\G^i(x(t_1)),\cdots,\pa_\nu U_\G^i(x(t_{2n-1}))\right)^\top,
\enn
and $W,\Psi,R_j(s),D_{\Gamma}$ are given in \eqref{0916-1} and \eqref{dgamma}.
As pointed out in Remark \ref{0916rem}, $W\Psi$ is viewed as an unknown vector to avoid the calculation of the inverse of $W$.
Multiplying $W$ on both sides of (\ref{0916-6}), we obtain an approximation of $2w'(s)|x'(w(s))|\pa_\nu U_{\Gamma,\pm}(x(w(s)))$ by
\be\label{240329-1}
D_{\Gamma}WV_\pm=(WH\mp D_{\Gamma})(W\Psi)+WD_{\Gamma}Y,
\en
where we have used the equality $WD_{\Gamma}=D_{\Gamma}W$.

Following \cite[Section 8.7]{Cakoni14}, we seek the solution to (\ref{0916-4})--(\ref{0916-5}) of the form
\be\label{240329-4}
v(x):=\int_{\G}\left\{\frac{\pa\Phi_k(x,y)}{\pa\nu(y)}[v](y)-\Phi_k(x,y)\left[\frac{\pa v}{\pa\nu}\right](y)\right\}ds(y),
\en
where the jumps are defined by $[v]:=v_+-v_-$ and $[\pa_\nu v]:=\pa_\nu v_+-\pa_\nu v_-$.
By jump relations and the boundary condition (\ref{0916-5}), we obtain the  boundary integral equations
\ben
[v]\!\!\!\!&=&\!\!\!\!-(\nu\cdot q)\pa_\nu U_{\G,+}+(\nu\cdot q)\pa_\nu U_{\G,-},\\
2\int_\G\Phi_k(x,y)\left[\frac{\pa v}{\pa\nu}\right](y)ds(y)\!\!\!\!&=&\!\!\!\!2\int_\G\frac{\pa\Phi_k(x,y)}{\pa\nu(y)}[v](y)ds(y)+[v]+2(\nu\cdot q)\pa_\nu U_{\G,+}(x),
\enn
which can be approximated by the linear systems
\be\label{240329-2}
X_1\!\!\!\!&=&\!\!\!\!Q(V_--V_+)/2,\\ \label{240329-3}
{\color{xxx}M_{corner}}WX_2\!\!\!\!&=&\!\!\!\!(LW+I)X_1+V_+,
\en
where
\ben
X_1\!\!\!\!&=&\!\!\!\!\left([v](x(t_0)),[v](x(t_1)),\cdots,[v](x(t_{2n-1}))\right)^\top,\\
X_2\!\!\!\!&=&\!\!\!\!\left([\pa_\nu v](x(t_0)),[\pa_\nu v](x(t_1)),\cdots,[\pa_\nu v](x(t_{2n-1}))\right)^\top,\\
Q\!\!\!\!&=&\!\!\!\!{\rm diag}((\nu\cdot q)(x(t_0)),(\nu\cdot q)(x(t_1)),\cdots,(\nu\cdot q)(x(t_{2n-1}))),\\
L\!\!\!\!&:=&\!\!\!\!\left(R_j(s_i)L_1(w(s_i),w(s_j))+\frac{\pi}nL_2(w(s_i),w(s_j))\right)_{i,j=0,1,\cdots,2n-1},\\
L_1(t,\tau)\!\!\!\!&=&\!\!\!\!\begin{cases}
\frac k{2\pi}\left\{x'_2(\tau)[x_1(\tau)-x_1(t)]-x'_1(\tau)[x_2(\tau)-x_2(t)]\right\}\frac{J_1^{(1)}(k|x(t)-x(\tau)|)}{|x(t)-x(\tau)|}, & t\neq\tau,\\
0, & t=\tau,
\end{cases}\\
L_2(t,\tau)\!\!\!\!&=&\!\!\!\!\begin{cases}
L(t,\tau)-L_1(t,\tau)\ln\left(4\sin^2\frac{t-\tau}2\right), & t\neq\tau,\\
-\frac1{2\pi}\frac{x''_2(t)x'_1(t)-x''_1(t)x'_2(t)}{|x'(t)|^2}, & t=\tau,\\
\end{cases}\\
L(t,\tau)\!\!\!\!&:=&\!\!\!\!\begin{cases}
-\frac{ik}2\left\{x'_2(\tau)[x_1(\tau)-x_1(t)]-x'_1(\tau)[x_2(\tau)-x_2(t)]\right\}\frac{H_1^{(1)}(k|x(t)-x(\tau)|)}{|x(t)-x(\tau)|}, & t\neq\tau,\\
-\frac1{2\pi}\frac{x''_2(t)x'_1(t)-x''_1(t)x'_2(t)}{|x'(t)|^2}, & t=\tau,
\end{cases}
\enn
and $V_\pm$ are given in (\ref{0916-6}), ${\color{xxx}M_{corner}},R_j(s),W$ are given in \eqref{0916-1}.
With the above notations in discrete form, we deduce from \eqref{240329-4} that
\ben
v(x)\!\!\!\!&\approx&\!\!\!\!\frac\pi n\left(\begin{array}{cccc}
\frac{\pa\Phi_k(x,x(t_0))}{\pa\nu(x(t_0))} & \frac{\pa\Phi_k(x,x(t_1))}{\pa\nu(x(t_1))} & \cdots & \frac{\pa\Phi_k(x,x(t_{2n-1}))}{\pa\nu(x(t_{2n-1}))}
\end{array}\right)D_{\Gamma}WX_1\\
&&-\frac\pi n\left(\begin{array}{cccc}
\Phi_k(x,x(t_0)) & \Phi_k(x,x(t_1))& \cdots & \Phi_k(x,x(t_{2n-1}))
\end{array}\right)D_{\Gamma}WX_2,\quad x\in\R^2\ba\ov{\G},\\
v^\infty(\hat x)\!\!\!\!&\approx&\!\!\!\!\frac\pi n\frac{e^{-i\frac\pi4}k}{\sqrt{8k\pi}}\left(\begin{array}{cccc}
e^{-ik\hat x\cdot x(t_0)}\hat x\cdot\nu(x(t_0)) & \cdots & e^{-ik\hat x\cdot x(t_{2n-1})}\hat x\cdot\nu(x(t_{2n-1}))
\end{array}\right)D_{\Gamma}WX_1\\
&&-\frac\pi n\frac{e^{i\frac\pi4}}{\sqrt{8k\pi}}\left(\begin{array}{cccc}
e^{-ik\hat x\cdot x(t_0)} & e^{-ik\hat x\cdot x(t_1)}& \cdots & e^{-ik\hat x\cdot x(t_{2n-1})}
\end{array}\right)D_{\Gamma}WX_2,\quad\hat x\in\Sp.
\enn
Analogously to Remark \ref{0916rem}, $WX_1$ and $WX_2$ can be viewed as two unknown vectors to avoid the calculation of the inverse of $W$.
Multiplying $W$ on both sides of \eqref{240329-2} and \eqref{240329-3} leads to
\be\no
WX_1\!\!\!\!&=&\!\!\!\!Q(WV_--WV_+)/2,\\ \label{240330-1}
W{\color{xxx}M_{corner}}(WX_2)\!\!\!\!&=&\!\!\!\!(WL+I)(WX_1)+WV_+,
\en
where $WV_\pm=D_{\Gamma}^{-1}[(WH\mp D_{\Gamma})(W\Psi)+WD_{\Gamma}Y]$ by \eqref{240329-1}.
Since it is not stable to solve \eqref{240330-1}, we may apply the Tikhonov regularization scheme to obtain
\be\label{211013-2}
WX_2\approx[\alpha_0I+(W{\color{xxx}M_{corner}})^*(W{\color{xxx}M_{corner}})]^{-1}(W{\color{xxx}M_{corner}})^*[(WL+I)(WX_1)+WV_+],
\en
or
\ben
WX_2\approx{\color{xxx}M_{corner}^{-1}}(\alpha_0I+W^*W)^{-1}W^*[(WL+I)(WX_1)+WV_+],
\enn
where the regularization parameter $\alpha_0\!>\!0$ is appropriately chosen.

Noting that this paper focuses on piecewise linear cracks, we update the locations of crack tips and interior corners $\{P_{\ell}:\ell=0,1,\cdots,N\}$ of the crack given in terms of (\ref{230912-1}) in each iteration step instead of the coefficients of basis shape functions such as Chebyshev polynomials in \cite{Kress1995} (see also \cite[Section 5.4]{CK19}).
{\color{xxx}Precisely, in the $m$-th iteration step the updated crack tips and interior corners $\{P_{\ell}^{(m)}=(P_{\ell,1}^{(m)},P_{\ell,2}^{(m)}):\ell=0,1,\cdots,N\}$ are given by
\be\label{240309-1}
\left(\begin{array}{c}
P_{0,1}^{(m)}\\
\vdots\\
P_{N,1}^{(m)}\\
P_{0,2}^{(m)}\\
\vdots\\
P_{N,2}^{(m)}
\end{array}\right)=\left(\begin{array}{c}
\widetilde P_{0,1}^{(m-1)}\\
\vdots\\
\widetilde P_{N,1}^{(m-1)}\\
\widetilde P_{0,2}^{(m-1)}\\
\vdots\\
\widetilde P_{N,2}^{(m-1)}
\end{array}\right)+\left(\begin{array}{c}
\Delta P_{0,1}^{(m)}\\
\vdots\\
\Delta P_{N,1}^{(m)}\\
\Delta P_{0,2}^{(m)}\\
\vdots\\
\Delta P_{N,2}^{(m)}
\end{array}\right),\quad m=1,2,\cdots,
\en
where $\{\widetilde P_{\ell}^{(m-1)}=(\widetilde P_{\ell,1}^{(m-1)},\widetilde P_{\ell,2}^{(m-1)}):\ell=0,1,\cdots,N\}$ are the crack tips and interior corners of $\G_{m-1}$, and $\{\Delta P_{\ell}^{(m)}=(\Delta P_{\ell,1}^{(m)},\Delta P_{\ell,2}^{(m)}):\ell=0,1,\cdots,N\}$ are given by \eqref{230912-1} and \eqref{211014-1}.
Here, $\G_0$ is the initial guess and $\G_{m-1}$ is the result of the $(m-1)$-th tangential update (see Subsection \ref{sub5.2} below).}

\subsection{Tangential updates of two crack tips}\label{sub5.2}

As shown by \cite[(45) and Theorem 5.2]{Bochniak2002}, the Frech\'et derivative at interior corners can be calculated via \eqref{0915}, while the Frech\'et derivative at two crack tips contains another terms related to the tangential updates.
To solve {\bf Problem 2}, we insert the tangential update after \eqref{240309-1} in each iteration step.

{\color{xxx}To introduce the tangential update, in the $m$-th iteration step with the notations in (\ref{240309-1}) define the following two tangential unit vectors:
\ben
\tau_0^{(m)}=\frac{P_1^{(m)}-P_0^{(m)}}{\left|P_1^{(m)}-P_0^{(m)}\right|},\quad\tau_N^{(m)}=\frac{P_N^{(m)}-P_{N-1}^{(m)}}{\left|P_N^{(m)}-P_{N-1}^{(m)}\right|}.
\enn
Set $P^{(m)}_{0,\pm}:=P_{0}^{(m)}\pm l_m\tau_0^{(m)}$ and $P^{(m)}_{N,\pm}:=P_N^{(m)}\pm l_m\tau_N^{(m)}$ with
\ben
l_m=\left(\frac1{N+1}\sum_{j=0}^N\left\{\left|\Delta P_{j,1}^{(m)}\right|^2+\left|\Delta P_{j,2}^{(m)}\right|^2\right\}\right)^{1/2}.
\enn
Let $\G_m^{(\pm,\pm)}$ denote the crack with tips (different from the crack given by (\ref{240309-1})) and interior corners (the same as the crack given by (\ref{240309-1})) given in order by $\{P^{(m)}_{0,\pm},P_1^{(m)},\cdots,P_{N-1}^{(m)},P^{(m)}_{N,{\pm}}\}$, respectively.
Let the crack $\G_m\in\{\G_m^{(\pm,\pm)}\}$ corresponding to the smallest residue among $\|U_{\G_m^{(\pm,\pm)}}^\infty-U_\G^\infty\|_{L^2(\Sp_0)}$ be the tangential update result of $m$-th iteration step given by (\ref{240309-1}).}

\section{Numerical examples}\label{sec4}
\setcounter{equation}{0}

In this section, we will display some numerical examples.
Let $m$ belong to a certain index set.
In the figures of numerical results to $I_{crack}$ defined by (\ref{comp0}), the colors for cracks  {\color{xxx}$\G_{a^{(m)}}$} are given in terms of RGB values in $[0,1]\times[0,1]\times[0,1]$:
\ben
\begin{cases}
(V^{(m)},1-V^{(m)},0) & \text{if}\;V^{(m)}\ge0,\\
(0,1+V^{(m)},-V^{(m)}) & \text{if}\;V^{(m)}<0,
\end{cases}
\quad\text{ with }
V^{(m)}:=2\times\frac{I_{crack}({\color{xxx}a^{(m)}})-I_{min}}{I_{max}-I_{min}}-1,
\enn
where $I_{max}:=\max_{m}\{I_{crack}({\color{xxx}a^{(m)}})\}$ and $I_{min}:=\min_{m}\{I_{crack}({\color{xxx}a^{(m)}})\}$.
In the figures of numerical results to (\ref{comparison1}), the colors for disks $B_R(P^{(m)})$ are given in terms of RGB values in a similar manner {\color{xxx}to} (\ref{comp0}).
{\color{xxx}However,} different from (\ref{comp0}) and (\ref{comparison1}), the numerical results to (\ref{comparison2}) are given by the Matlab colormap 'jet'.
For the numerical results of the factorization method with a single wave, in the figures of numerical results to (\ref{I-1}), (\ref{211005-3}), (\ref{al3.1}) and (\ref{211012-2}), the colors for test domains $D^{(m)}$ (with different boundary conditions or refractive indices) are also given in terms of RGB values in a similar manner {\color{xxx}to $I_{crack}$ defined by} (\ref{comp0}).
Furthermore, for Newton's iteration method, the color of crack $\G_m$ in $m$-th iteration step is given in terms of RGB values in $[0,1]\times[0,1]\times[0,1]$:
\ben
\begin{cases}
(V^{(m)},1-V^{(m)},0) & \text{if}\;V^{(m)}\ge0,\\
(0,1+V^{(m)},-V^{(m)}) & \text{if}\;V^{(m)}<0,
\end{cases}
\quad\text{ with }
V^{(m)}:=\frac{2m}{\text{Total iteration number}}-1.
\enn
\begin{example}\label{example0} \textbf{(Numerical results of (\ref{comp0}))}
Let $\G$ be a sound-soft piecewise linear crack with two tips and an interior corner given in order by $(1,3)$, $(3,1)$ and $(2,0)$.
{\color{xxx}The measured data are the limited aperture far-field patterns $\{w^\infty(\hat x_p,y_0;k_0):p=0,1,\cdots,L\}$ with $k_0=5$,} $y_0=(8,0),(0,8),(-8,0),(0,-8)$, $\hat x_p=(\cos\theta_p,\sin\theta_p)$, $\theta_p=\frac\pi 2+\frac{\pi}Lp$, and $L=128$.
The numerical results of (\ref{comp0}) with {\color{xxx}$a=a^{(m_1,m_2)}=(2.3m_1-0.1,2.3m_2-0.1)$, $m_1,m_2\in\{-3,-2,-1,0,1\}$,} are shown in Figure \ref{fig240322}.

\begin{figure}[htb]
  \centering
  \subfigure[$y_0=(8,0)$]{
  \includegraphics[width=0.23\textwidth]{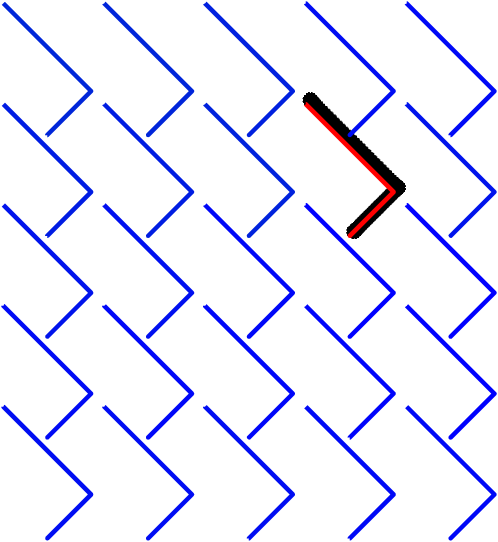}}
  \subfigure[$y_0=(0,8)$]{
  \includegraphics[width=0.23\textwidth]{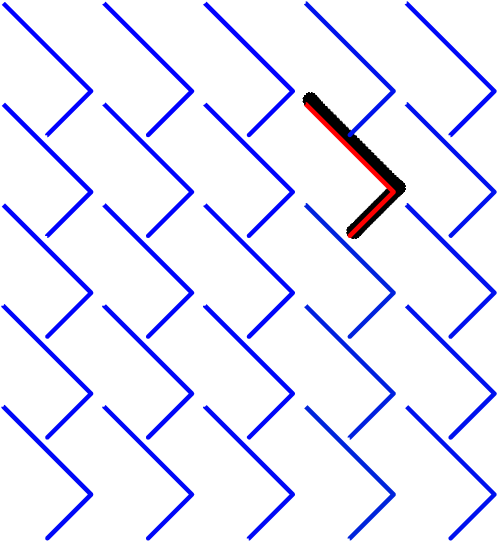}}
  \subfigure[$y_0=(-8,0)$]{
  \includegraphics[width=0.23\textwidth]{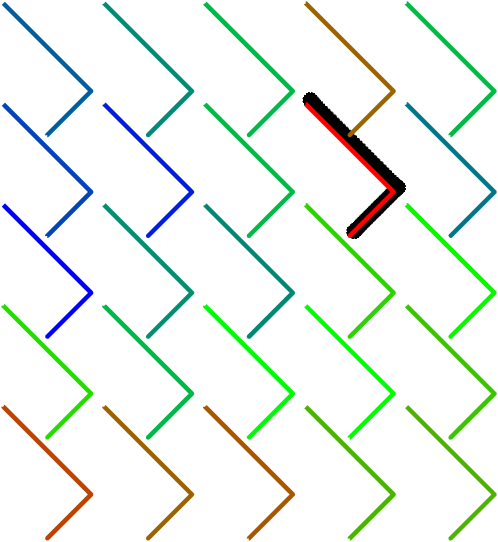}}
  \subfigure[$y_0=(0,-8)$]{
  \includegraphics[width=0.23\textwidth]{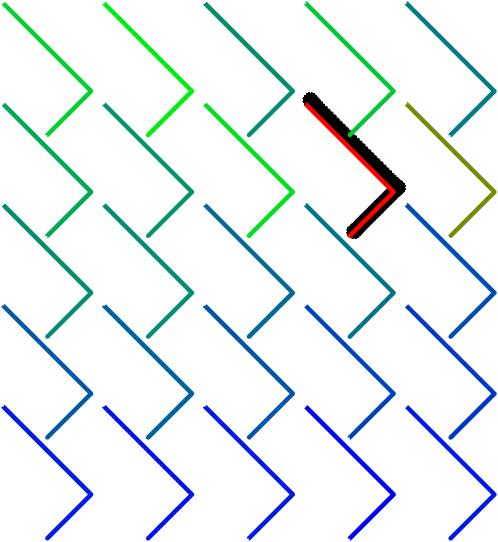}}
  \caption{Numerical results for Example \ref{example0} with $\G$ denoted by the black thick line.}\label{fig240322}
\end{figure}
\end{example}

\begin{example}\label{example1}
Let $\G$ be a sound-soft piecewise linear crack with two tips and an interior corner given in order by $(1,3)$, $(3,1)$ and $(2,0)$.

\textbf{(a) ({\color{xxx}Determine a rough location} by (\ref{comparison2})).}
{\color{xxx}The measured data are $\{u^\infty(\hat x_p,d_0;k_0):p=0,1,\cdots,L-1\}$ with $k_0=1,2,4,8$,} $d_0=(1,0)$, $\hat x_p=(\cos\theta_p,\sin\theta_p)$, $\theta_p=\frac{2\pi}Lp$, and $L=128$.
The numerical results for (\ref{comparison2}) with $\tau=1$ are shown in Figure \ref{fig230802-2}.

\begin{figure}[htb]
  \centering
  \subfigure[$k_0=1$]{
  \includegraphics[width=0.23\textwidth]{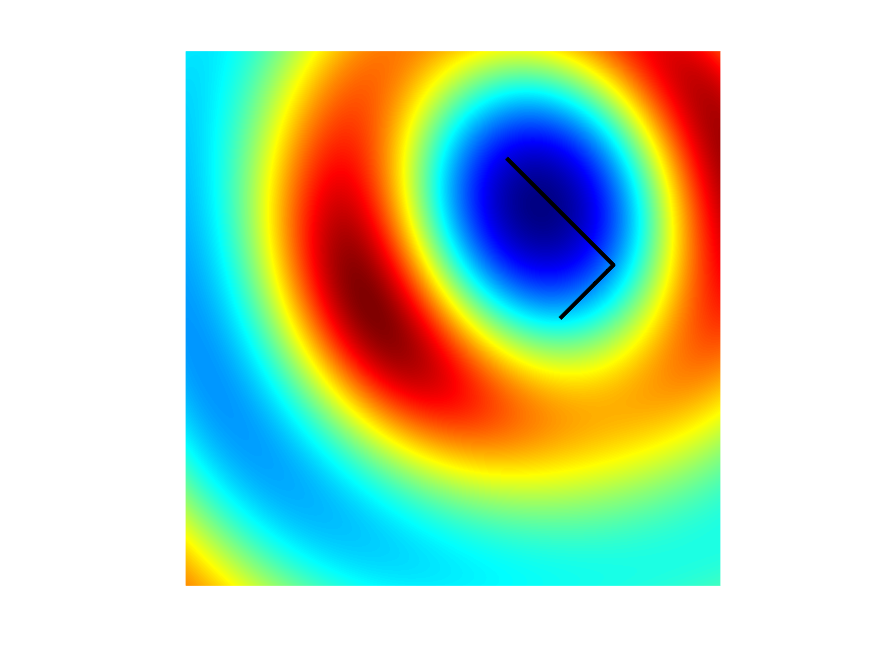}}
  \subfigure[$k_0=2$]{
  \includegraphics[width=0.23\textwidth]{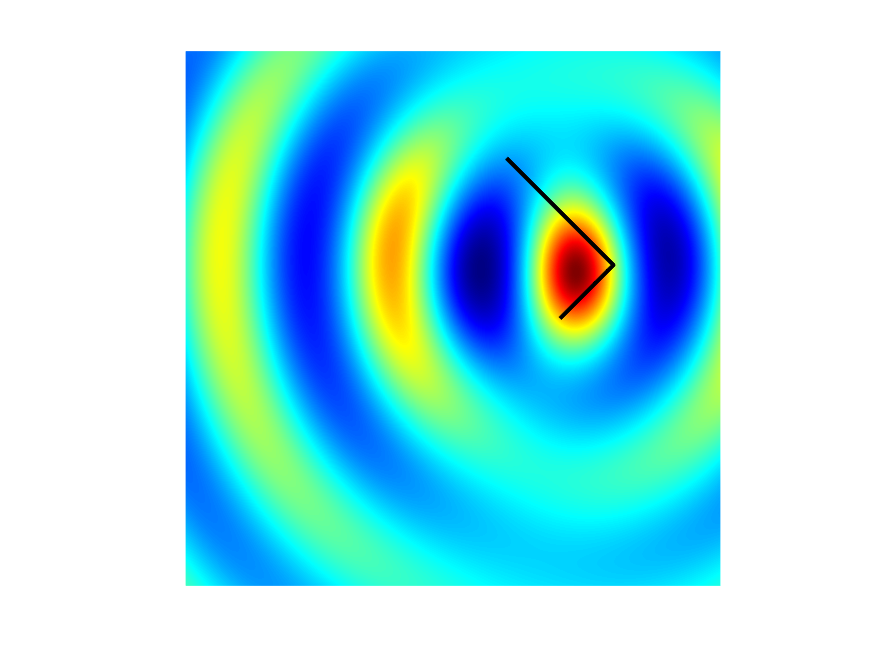}}
  \subfigure[$k_0=4$]{
  \includegraphics[width=0.23\textwidth]{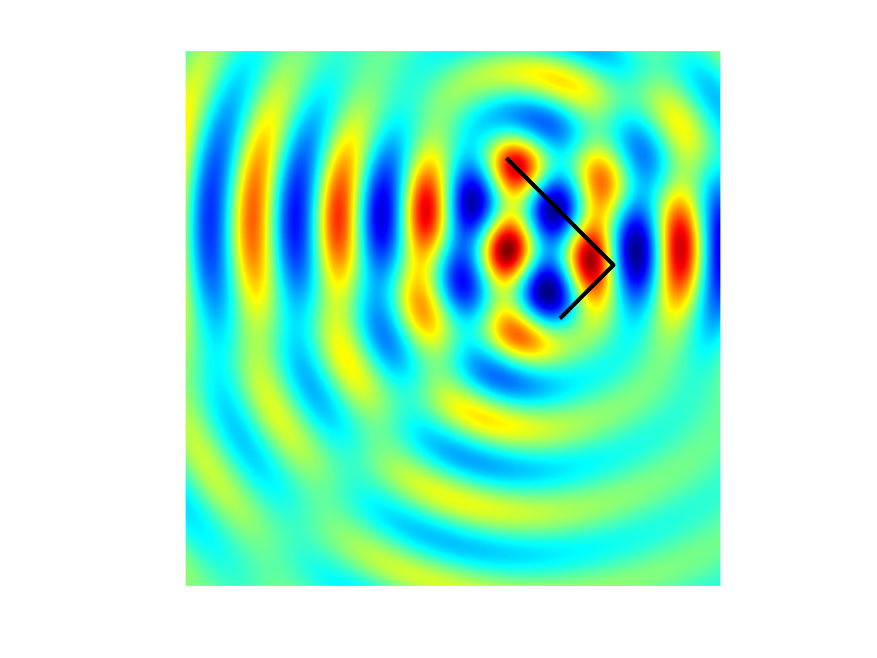}}
  \subfigure[$k_0=8$]{
  \includegraphics[width=0.23\textwidth]{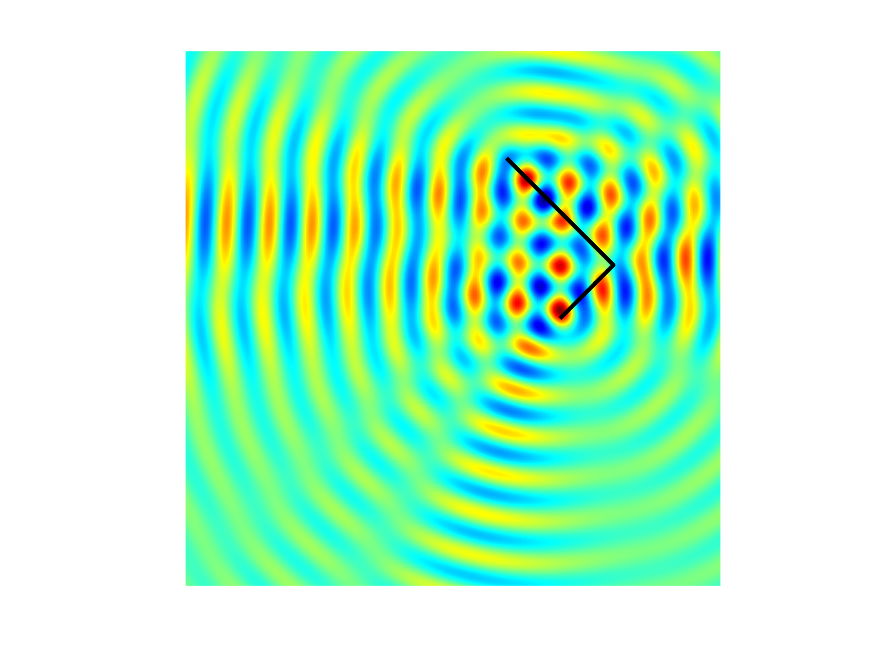}}
  \caption{Numerical results for Example \ref{example1} (a) with $\G$ denoted by the black solid line.}\label{fig230802-2}
\end{figure}

\textbf{(b) ({\color{xxx}Determine a rough location} by (\ref{comparison1})).}
{\color{xxx}The measured data are $\{u^\infty(\hat x_p,d_0;k_0):p=0,1,\cdots,L-1\}$ with $k_0=2$,} $d_0=(1,0),(0,1),(-1,0),(0,-1)$, $\hat x_p=(\cos\theta_p,\sin\theta_p)$, $\theta_p=\frac{2\pi}Lp$, and $L=128$.
The numerical results for (\ref{comparison1}) with sound-soft disks of {\color{xxx} the same radius $0.2$ and different centers} are shown in Figure \ref{fig230802-1}.

\begin{figure}[htb]
  \centering
  \subfigure[$d_0=(1,0)$]{
  \includegraphics[width=0.23\textwidth]{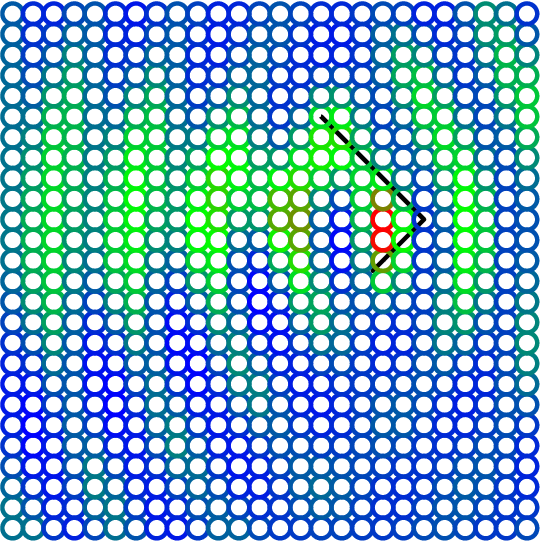}}
  \subfigure[$d_0=(0,1)$]{
  \includegraphics[width=0.23\textwidth]{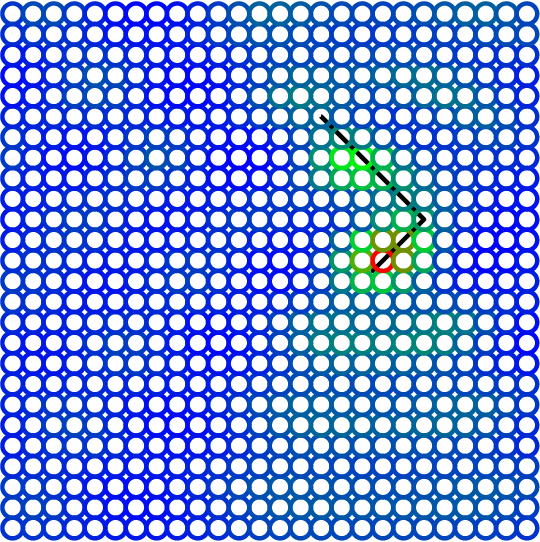}}
  \subfigure[$d_0=(-1,0)$]{
  \includegraphics[width=0.23\textwidth]{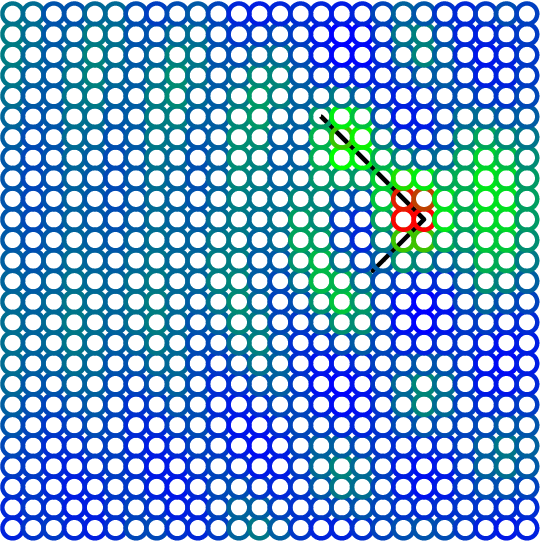}}
  \subfigure[$d_0=(0,-1)$]{
  \includegraphics[width=0.23\textwidth]{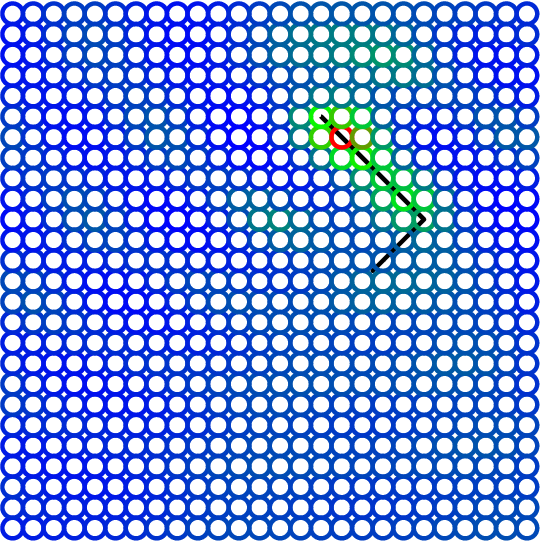}}
  \caption{Numerical results for Example \ref{example1} (b) with $\G$ denoted by the black solid line.}\label{fig230802-1}
\end{figure}
\end{example}

\begin{example}\label{example2}
Let $\G$ be a sound-soft piecewise linear crack with two tips and interior corners given in order by $(0,2)$, $(-1,1)$, $(1,-1)$, $(0,-2)$.

\textbf{(a) (Factorization method with a single wave for test {\color{xxx}scatterers} of different boundary conditions and different refractive indices).}
{\color{xxx}The measured data are} $\{u^\infty(\hat x_p,d_0;k_0):p=0,1,\cdots,L-1\}$ with $k_0=2$, $d_0=(1,0)$, $\hat x_p=(\cos\theta_p,\sin\theta_p)$, $\theta_p=\frac{2\pi}Lp$, and $L=64$.
For each $m\in\{5,\cdots,25\}$, set the test domain $D=D^{(m)}_{D}$, $D^{(m)}_{N}$, $D^{(m)}_{I}$, $D^{(m)}_{n}$ to be a sound-soft disk, a sound-hard disk, an impedance disk with impedance coefficient $\eta\!=\!ik_0$, a medium with constant refractive index $m$ in the disk, respectively, centered at $P\!=\!(0,1)$ with radius $m/5$.
The numerical results of (\ref{I-1}) with test domains given as above are shown in Figure \ref{fig8}.
\begin{figure}[htb]
  \centering
  \subfigure[$I_1(D^{(m)}_{D})$]{
  \includegraphics[width=0.2\textwidth]{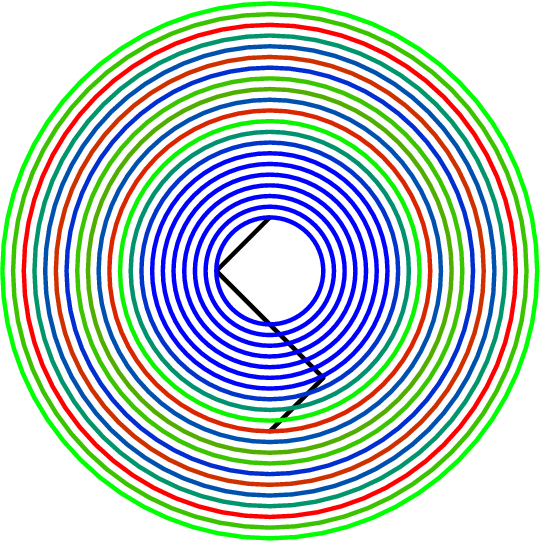}}
  \subfigure[$I_1(D^{(m)}_{N})$]{
  \includegraphics[width=0.2\textwidth]{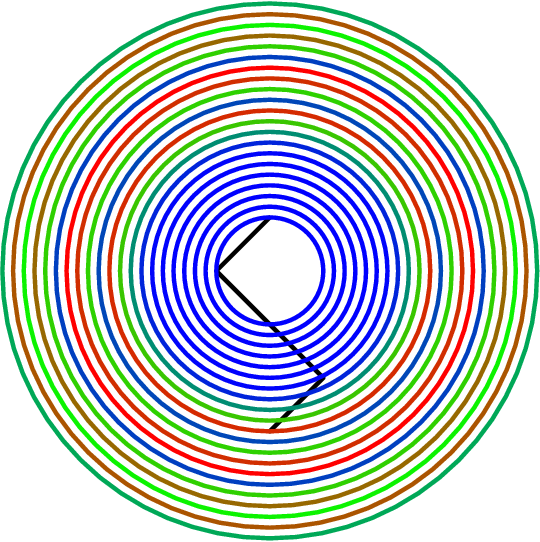}}
  \subfigure[$I_1(D^{(m)}_{I})$]{
  \includegraphics[width=0.2\textwidth]{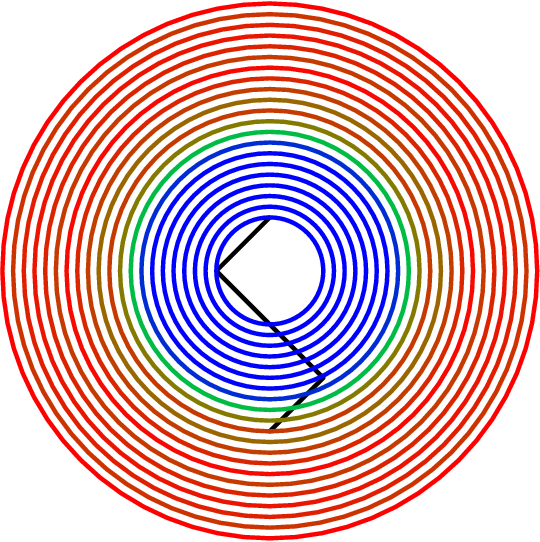}}\\
  \subfigure[$I_1(D^{(m)}_{n=1\!/\!4})$]{
  \includegraphics[width=0.2\textwidth]{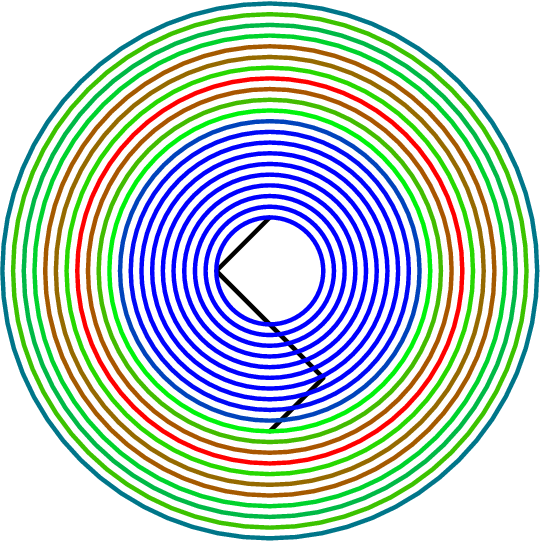}}
  \subfigure[$I_1(D^{(m)}_{n=4})$]{
  \includegraphics[width=0.2\textwidth]{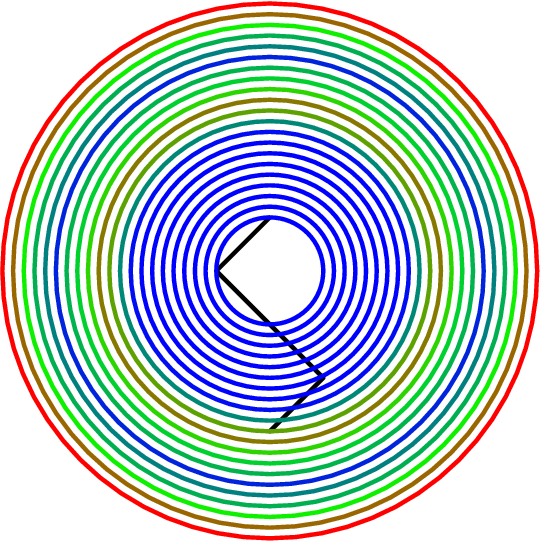}}
  \subfigure[$I_1(D^{(m)}_{n=3\!+\!4i})$]{
  \includegraphics[width=0.2\textwidth]{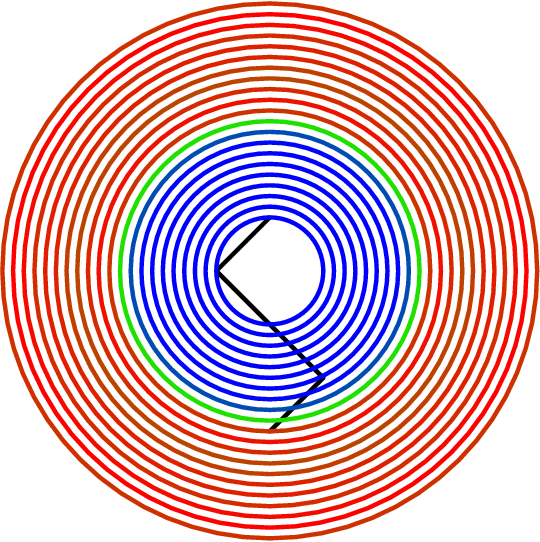}}
  \caption{Numerical results for Example \ref{example2} (a) with $\G$ denoted by the black solid line.}\label{fig8}
\end{figure}

\textbf{(b) (Point source incidence).}
{\color{xxx}The measured data are} $\{w^\infty(\hat x_p,y_0;k_0):p=0,1,\cdots,L\}$ with $k_0=2$, $y_0=(2,-1)$, $\hat x_p=(\cos\theta_p,\sin\theta_p)$, $\theta_p=\frac{2\pi}Lp$, and $L=64$.
For each $m\in\{5,\cdots,25\}$, set the test domain $D=D^{(m)}_P$ to be an impedance disk centered at $P$ of radius $m/5$ with impedance coefficient $\eta=ik_0$.
The numerical results of (\ref{I-1}) with $D=D^{(m)}_P$ are given by Figure \ref{fig0918}.

\begin{figure}[htb]
  \centering
  \subfigure[$D_{P=(-3,0)}^{(m)}$]{
  \includegraphics[width=0.2\textwidth]{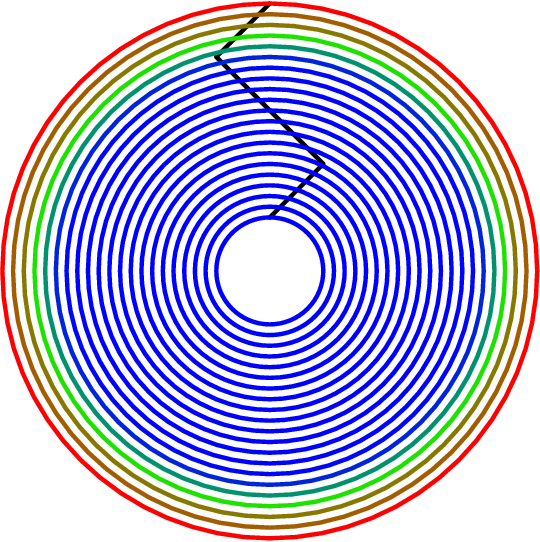}}
  \subfigure[$D_{P=(-1,0)}^{(m)}$]{
  \includegraphics[width=0.2\textwidth]{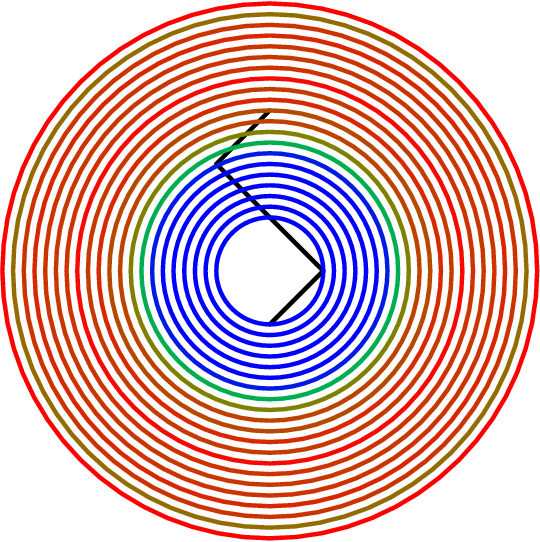}}
  \subfigure[$D_{P=(1,0)}^{(m)}$]{
  \includegraphics[width=0.2\textwidth]{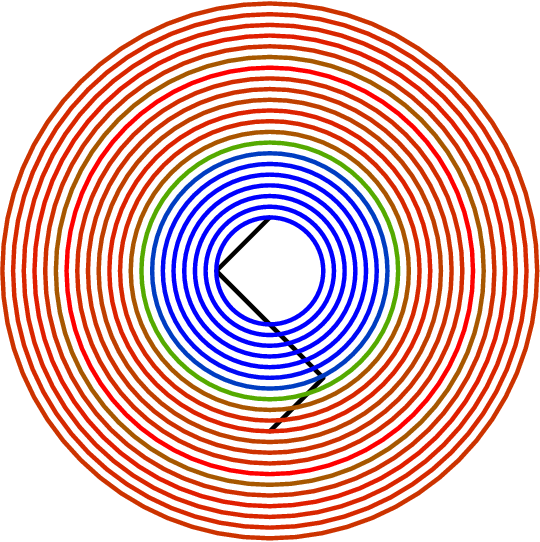}}
  \subfigure[$D_{P=(3,0)}^{(m)}$]{
  \includegraphics[width=0.2\textwidth]{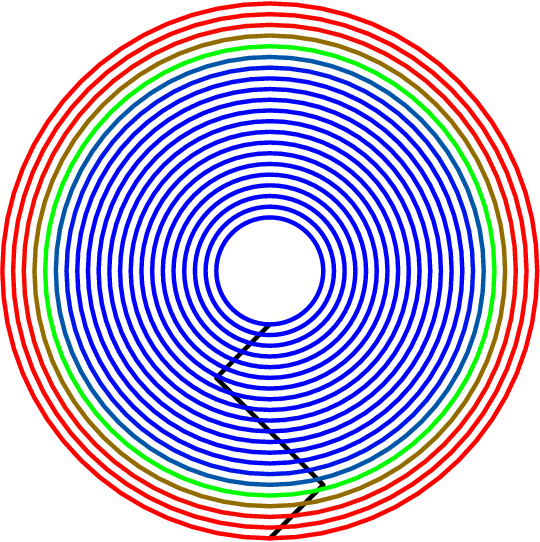}}
  \caption{Numerical results for Example \ref{example2} (b) with $\G$ denoted by the black solid line.}\label{fig0918}
\end{figure}
\end{example}
{\color{xxx}We remind the reader that $k^2>0$ is not an impedance eigenvalue of $-\Delta$ in $D$ provided $D$ is an impedance obstacle with ${\rm Im}\,\eta\neq0$ on $\pa D$ (see \cite[Theorem 8.2]{Cakoni14}) and $k^2>0$ is not an interior transmission eigenvalue in $D$ provided ${\rm Im}\,n\neq0$ and $\inf_{x\in\ov{D}}|n(x)-1|>0$ in the compact support $D$ of $n-1$ (see \cite[Theorem 8.12]{CK19}).
In view of Theorems \ref{thm3.9} and \ref{thm3.9-la} and numerical results in Figures \ref{fig8} and \ref{fig0918}, we always set the test scatterer $D$ to be an impedance obstacle or an inhomogeneous medium with a refractive index of nonzero imaginary part in the one-wave factorization method, to avoid the influence of possible eigenvalues.}

We are now ready to introduce our {\color{xxx}\textbf{hybrid method}} to detect a piecewise linear crack, which can be divided into the following three steps:

\textbf{Step 1: Detection for a rough location.}
We have the following two different approaches:

\textbf{Approach 1. (Contrast sampling method).} The rough location of target crack can be recovered by either (\ref{comparison2}) or (\ref{comparison1}) as shown in Figures \ref{fig230802-2} and \ref{fig230802-1}.

\textbf{Approach 2.}
Let $P_1,P_2,P_3$ be three different points such that they are not located on the same straight line.
Set an appropriate step length $r>0$ and {\color{xxx}a sufficiently large number} $N\in\Z_+$.
For each $m\in\{1,\cdots,N\}$, set $D_P^{(m)}$ to be a disk centered at $P\in\{P_1,P_2,P_3\}$ with radius $m r$.
Calculate the values of indicator functions (\ref{I-1}) for $D_P^{(m)}$.
As explained in Section \ref{sec3} and shown by Figures \ref{fig8} and \ref{fig0918}, the target crack is located near the critical circle whose neighbouring circles of the same center admit a big change of values to corresponding indicator functions.
Therefore, the rough location of the unknown crack is given by the intersection of the three critical circles.
This approach is verified by Example \ref{example3} (a) below.

\textbf{Step 2: {\color{xxx}Detection for a precise location and the convex hull of the crack.}}
For each point $P\in\R^2$, set $r_P:=\max\{r>0:{\color{xxx}I_1(B_r(P))}=0\}$ with $I_1$ defined by (\ref{I-1}), where {\color{xxx}$B_r(P)$} is the test disk centered at $P$ with radius $r>0$.
Furthermore, we define $\chi_P(x)=1$ if $|x-P|\leq r_P$ and $\chi_P(x)=0$ if $|x-P|>r_P$.
Theoretically, the convex hull of the crack is thus contained in the set of the maximum points of $\sum_{j\in J}\chi_{P_j}(x)$, where $\{P_j:j\in J\}$ are points located in $\R^2$.
It can be easily seen that the convex hull of the crack can be approximated by the set of the maximum points of $\sum_{j\in J}\chi_{P_j}(x)$ provided there are sufficiently many points $\{P_j:j\in J\}$.
The numerical results of this step are shown in Example \ref{example3} (b) below (see also {\color{xxx}Example} \ref{example3+} for square-shaped test domains).

\textbf{Step 3: Iteration method for a precise shape.}
With an appropriate initial guess from Step 2, we can apply the iteration method introduced in Section \ref{iteration_method} for a precise shape reconstruction (see Example \ref{example3} (c) below).

If the full aperture data is replaced by limited aperture data, the above hybrid method also works with the indicator function (\ref{I-1}) replaced by (\ref{211005-3}) (see Example \ref{ex-la} below).

\begin{example}\label{example3}
Let $\G$ be a sound-soft piecewise linear crack with two tips and an interior corner given in order by $(-1,1)$, $(1,-1)$ and $(0,-2)$.
{\color{xxx}The far-field data are} $\{u^\infty(\hat x_p,d_0;k_0):p=0,1,\cdots,L-1\}$ with $k_0=5$, $d_0=(-1,0)$, and $\hat x_p=(\cos\theta_p,\sin\theta_p)$, $\theta_p=\frac{2\pi}Lp$, $L=800$.

\textbf{(a) Step 1 (Detection for a rough location).}
We will compare the numerical results for shifted cracks {\color{xxx}$\G_a:=\{x+a\in\R^2:x\in\G\}$ for different $a\in\R^2$}.
For each $m\in\{1,\cdots,15\}$, set $D^{(m)}_P$ to be an impedance disk with impedance coefficient $\eta=ik_0$ centered at $P$ of radius $m$.
The numerical results of (\ref{al3.1}) with $D=D^{(m)}_P$ and $\alpha=10^{-8}$ are given by Figure \ref{fig210820-2}, {\color{xxx}where numerical results for $P\in\{(-10,0),(10,0)\}$ are given in Figure \ref{fig210820-2} (a)--(e), and numerical results for $P\in\{10(\cos\theta,\sin\theta):\theta=-\frac\pi6,\frac\pi2,\frac{7\pi}6\}$ are given in Figure \ref{fig210820-2} (f)--(j).}
\begin{figure}[htb]
  \centering
  \subfigure[$a=(-3,0)$]{
  \includegraphics[width=0.18\textwidth]{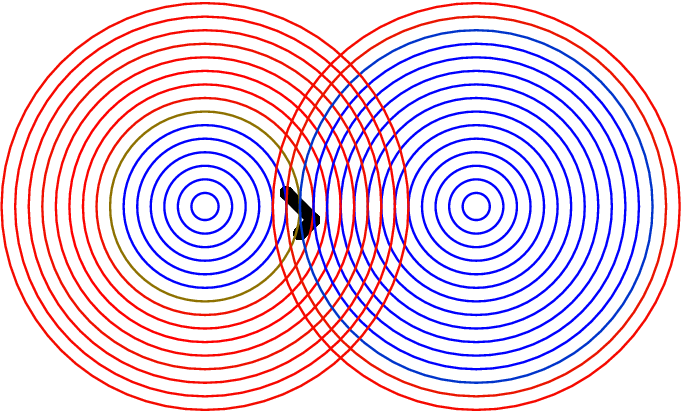}}
  \subfigure[$a=(-1,0)$]{
  \includegraphics[width=0.18\textwidth]{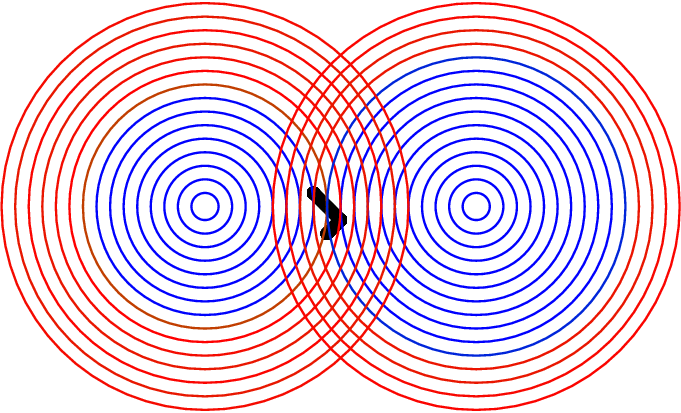}}
  \subfigure[$a=(1,0)$]{
  \includegraphics[width=0.18\textwidth]{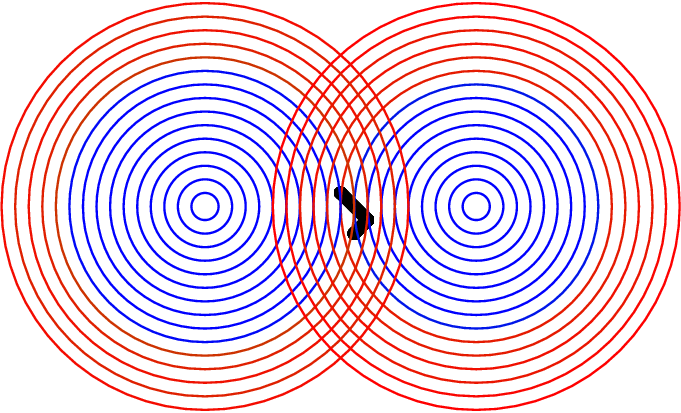}}
  \subfigure[$a=(3,0)$]{
  \includegraphics[width=0.18\textwidth]{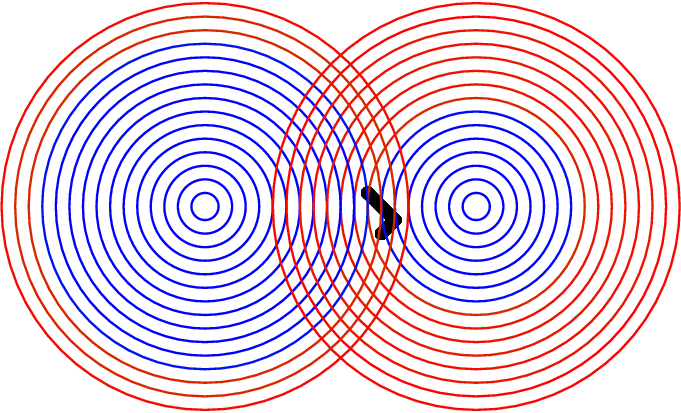}}
  \subfigure[$a=(5,0)$]{
  \includegraphics[width=0.18\textwidth]{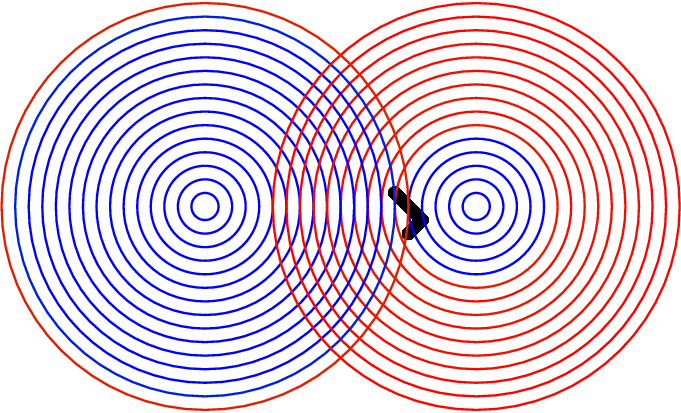}}
  \subfigure[$a=(-3,0)$]{
  \includegraphics[width=0.18\textwidth]{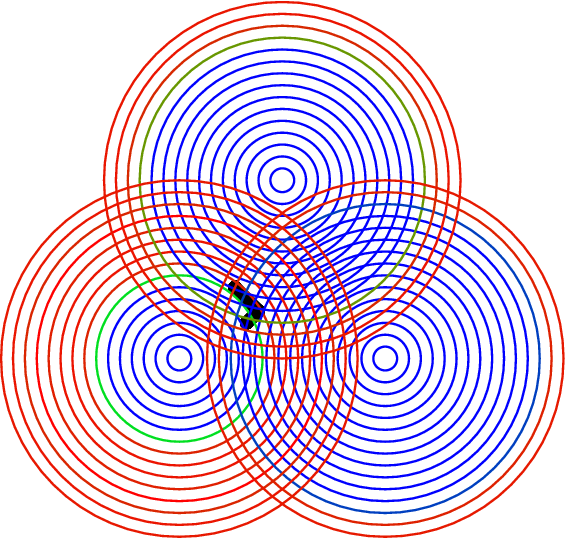}}
  \subfigure[$a=(-1,0)$]{
  \includegraphics[width=0.18\textwidth]{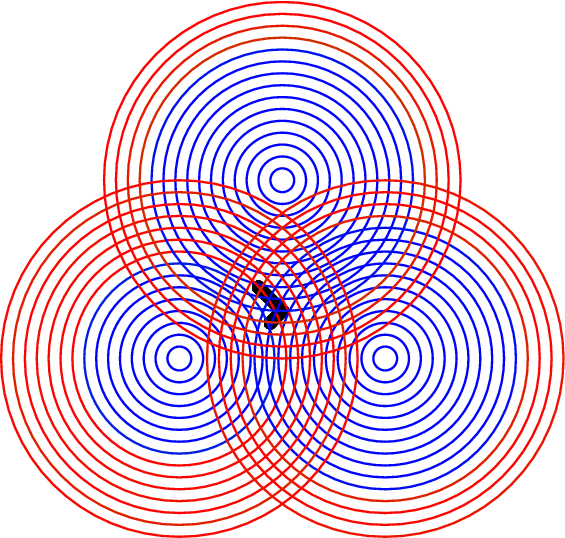}}
  \subfigure[$a=(1,0)$]{
  \includegraphics[width=0.18\textwidth]{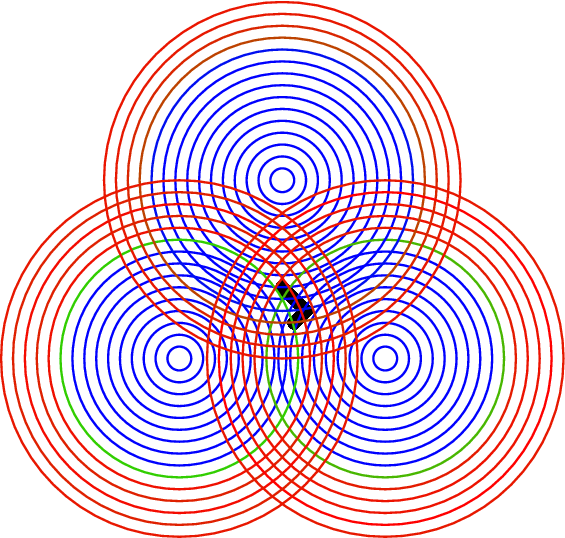}}
  \subfigure[$a=(3,0)$]{
  \includegraphics[width=0.18\textwidth]{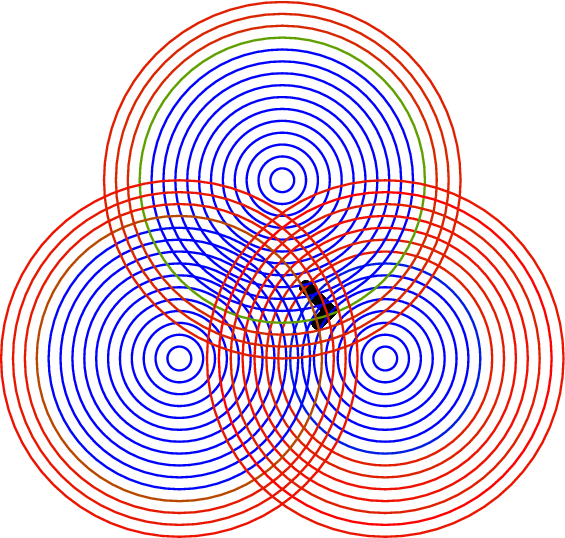}}
  \subfigure[$a=(5,0)$]{
  \includegraphics[width=0.18\textwidth]{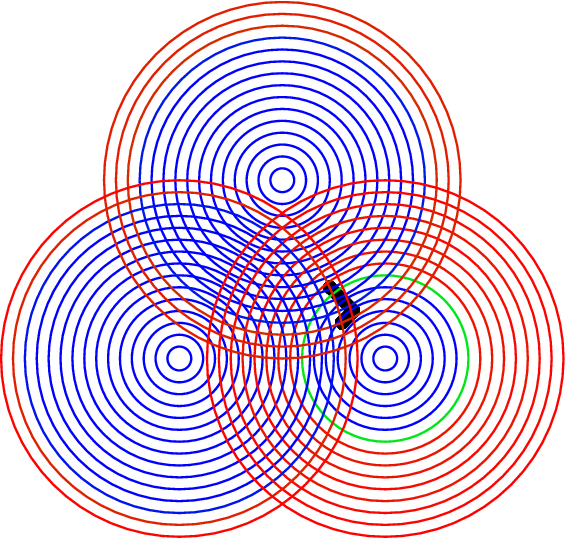}}
  \caption{Numerical results for Example \ref{example3} (a) {\color{xxx}with
   $\G_a$ denoted by the black thick line}.}\label{fig210820-2}
\end{figure}

\textbf{(b) Step 2 (Detection for a precise location and the convex hull of the crack).}
Following the strategy in Step 2 {\color{xxx}of the \textbf{hybrid method}}, we firstly set $P=(10,0)$ and calculate the values of {\color{xxx}$\widetilde{I}_1(B_r(P))$} defined by \eqref{al3.1} with $r=m/2$ for $m=2,3,\cdots,30$, where {\color{xxx}$B_r(P)$} denotes an impedance disk centered at $P$ with radius $r$ and constant impedance coefficient $\eta=ik_0$ (see Figures \ref{fig230803-1} (a) and (b)), where the regularization {\color{xxx}parameter in} (\ref{al3.1}) is {\color{xxx}set to be} $\alpha=10^{-8}$.
According to Figure \ref{fig230803-1} (b), we set the threshold to be $\epsilon=1.1\times10^{-5}$.
Secondly, we set $P_j=10(\cos\theta_j,\sin\theta_j)$ with $\theta_j=\frac{\pi}{16}j$, $j=0,1,\cdots,31$.
Define $m_j:=\max\{m\in\Z_+:{\color{xxx}\widetilde I_1(B_{0.1m}(P_j))}<\epsilon\}$, where the regularization {\color{xxx}parameter in} (\ref{al3.1}) is also {\color{xxx}set to be} $\alpha=10^{-8}$.
As an approximation of the function $\chi_P(x)$ in Step 2 {\color{xxx}of the \textbf{hybrid method}}, we define $\tilde\chi_{P_j}(x)=1$ if $|x-P_j|\leq 0.1m_j$ and $\tilde\chi_{P_j}(x)=0$ if $|x-P_j|> 0.1m_j$ for each $j$.
{\color{xxx}The numerical result for $\sum_{j=0}^{31}\tilde\chi_{P_j}(x)$ is} shown in Figure \ref{fig230803-1} (c), {\color{xxx}where the colors are given by the Matlab colormap 'jet'}.

\begin{figure}[htb]
  \centering
  \subfigure[]{
  \includegraphics[width=0.23\textwidth]{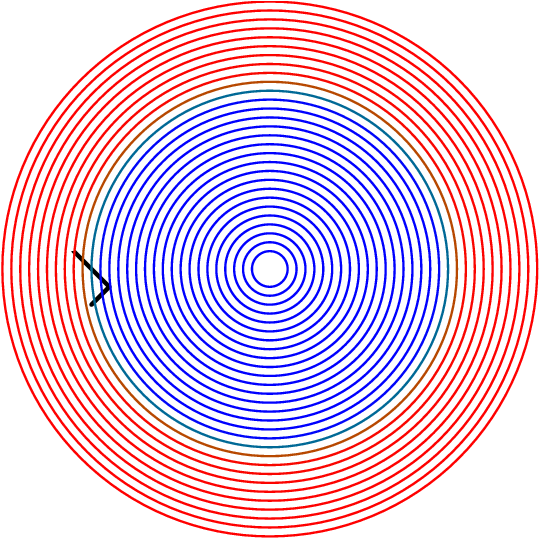}}\quad
  \subfigure[]{
  \includegraphics[width=0.3\textwidth]{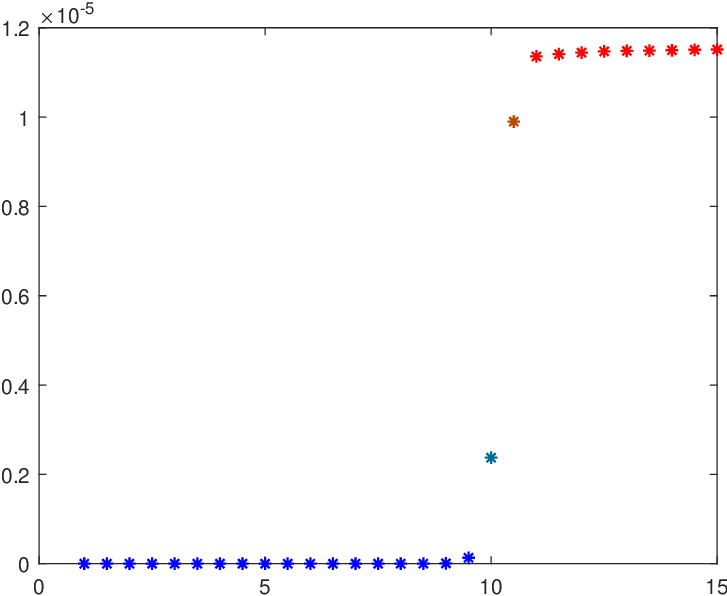}}\quad
  \subfigure[]{
  \includegraphics[width=0.35\textwidth]{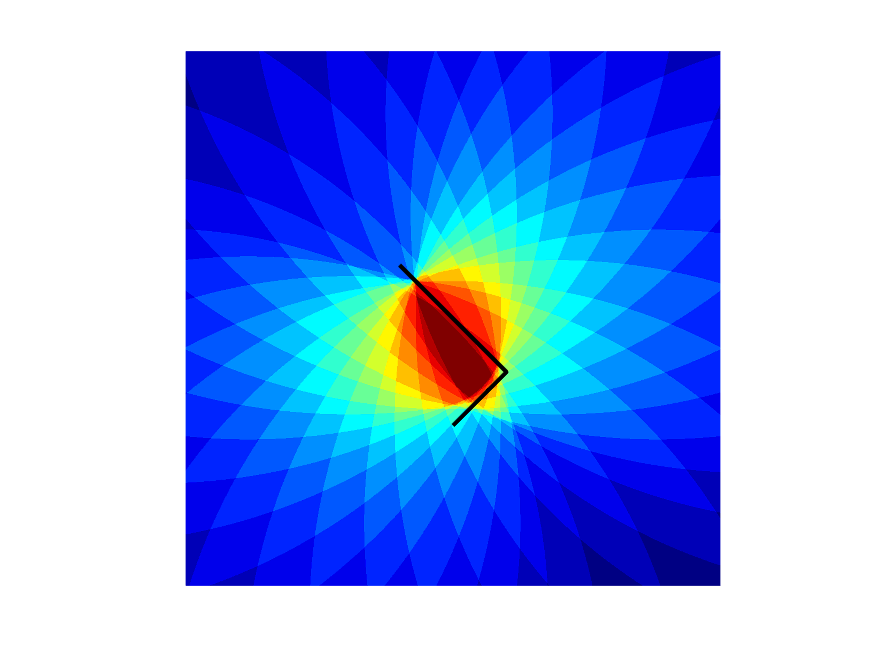}}
  \caption{Numerical results for Example \ref{example3} (b) {\color{xxx}with $\G$ denoted by the black solid line}.}\label{fig230803-1}
\end{figure}

\textbf{(c) Step 3 (Iteration for a precise shape).}
For a more precise numerical result, we apply the iteration method based on noisy far-field data $u_\delta^\infty(\hat x_p,d_0;k_0)$ given similarly to (\ref{noi}).
The initial guess is given by the location of crack tips and interior corners in order by $(-0.8,1.2)$, $(0.5,-1.2)$, $(-0.2,-1.9)$.
We choose different noise ratios $\delta\!>\!0$ as shown in Figure \ref{fig10} (a)--(c).
Noting that {\color{xxx}one cannot} take it for grant that $\G$ consists of two straight lines, we also give numerical results with initial guess given by $(-0.8,1.2)$, $(-0.2,-0.2)$, $(0.7,-1.2)$, $(-0.2,-1.9)$ as shown in Figure \ref{fig10} (d) and by $(-0.8,1.2)$, $(-0.2,-0.2)$, $(0.8,-0.8)$, $(0.3,-1)$, $(-0.2,-1.9)$ as shown in Figure \ref{fig10} (e), respectively.
The number of total iteration steps for each figure is $10$, and we set $\alpha=10$ in (\ref{211014-1}) and $\alpha_0=10^{-2}$ in (\ref{211013-2}).

\begin{figure}[htb]
  \centering
  \subfigure[$\delta=0\%$]{
  \includegraphics[width=0.18\textwidth]{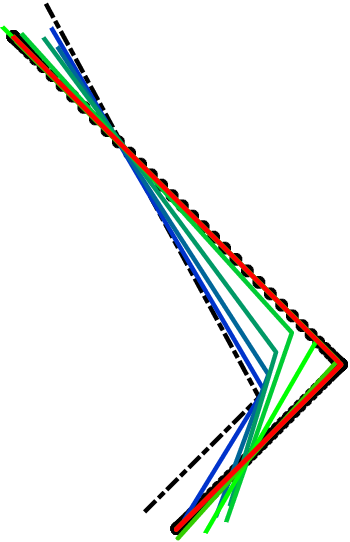}}
  \subfigure[$\delta=1\%$]{
  \includegraphics[width=0.18\textwidth]{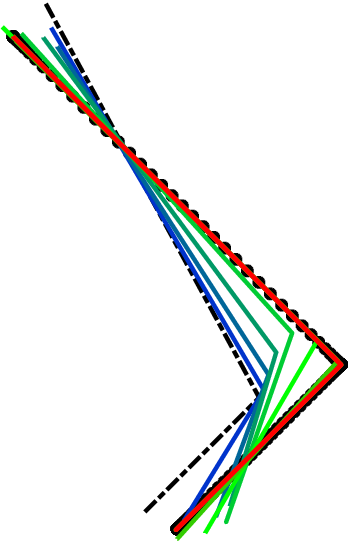}}
  \subfigure[$\delta=10\%$]{
  \includegraphics[width=0.18\textwidth]{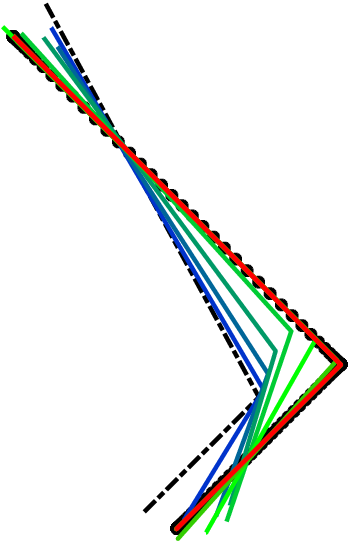}}
  \subfigure[$\delta=0\%$]{
  \includegraphics[width=0.18\textwidth]{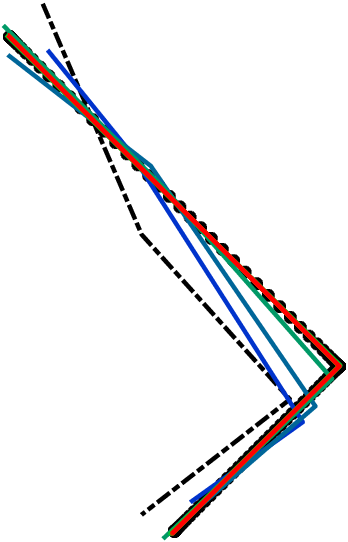}}
  \subfigure[$\delta=0\%$]{
  \includegraphics[width=0.18\textwidth]{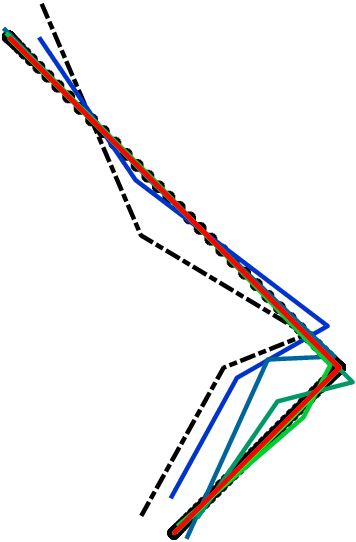}}
  \caption{{\color{xxx}Numerical results for} Example \ref{example3} (c). The black dots '*' represent the true crack {\color{xxx}$\G$}, the dashed line '-.' represents the initial guess, and colored lines represent the numerical result in each iteration step.}\label{fig10}
\end{figure}
\end{example}

It should be remarked that the numerical result of iteration method depends heavily on the initial guess.
If the initial guess is not sufficiently close to the true shape, then the iteration method may not give a satisfactory result.
{\color{xxx}As an advantage of our hybrid merthod,} the factorization method with a single wave can provide a quite good initial guess.

\begin{example}\label{example3+}
\textbf{(Square-shaped test domains).}
Let $\G$ be a sound-soft piecewise linear crack with two tips and an interior corner given in order by $(-1,1)$, $(1,-1)$ and $(0,-2)$.
{\color{xxx}The far-field data are} $\{u^\infty(\hat x_p,d_0;k_0):p=0,1,\cdots,L-1\}$ with {\color{xxx}$k_0=\frac14,\frac12,1,2,4$}, $d_0=(-1,0)$, and $\hat x_p=(\cos\theta_p,\sin\theta_p)$, $\theta_p=\frac{2\pi}Lp$, $L\!=\!64$.
For each $\ell\in\{2,\cdots,25\}$, set $D^{(m)}_P$ to be an impedance square with $\eta=ik_0$ centered at {\color{xxx} $P$} with side length $2m/5$.
The numerical results of (\ref{I-1}) {\color{xxx}with $D=D^{(m)}_P$ for different $P\in\R^2$} are given by Figure \ref{fig9}, where we choose different centers $P\in\R^2$ and different wave numbers {\color{xxx}$k_0>0$}.

\begin{remark}
With properties of the boundary integral operator {\color{xxx}for} Lipschitz domains (see \cite{Costabel}), the extension of the factorization method from domains of class $C^2$ to Lipschitz domains can be established similarly.
We refer the reader to \cite[Chapter 5]{KH2015} and \cite{M} for more
details on boundary value problems in Lipschitz domains. \end{remark}

\begin{figure}[htb]
  \centering
  \subfigure[$k_0\!=\!\frac14,D_{P=(3,0)}^{(m)}$]{
  \includegraphics[width=0.18\textwidth]{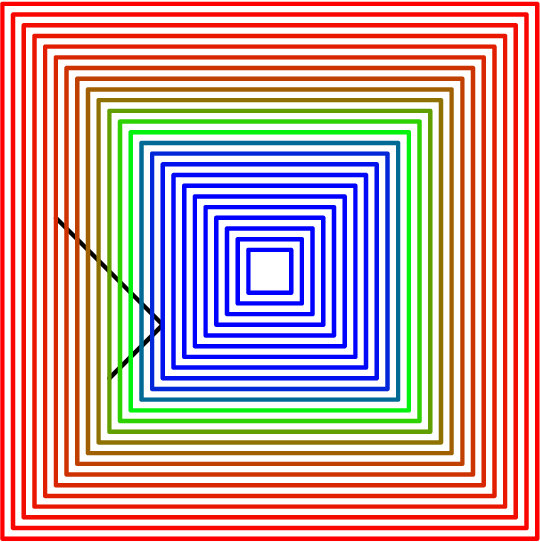}}
  \subfigure[$k_0\!=\!\frac12,D_{P=(2,0)}^{(m)}$]{
  \includegraphics[width=0.18\textwidth]{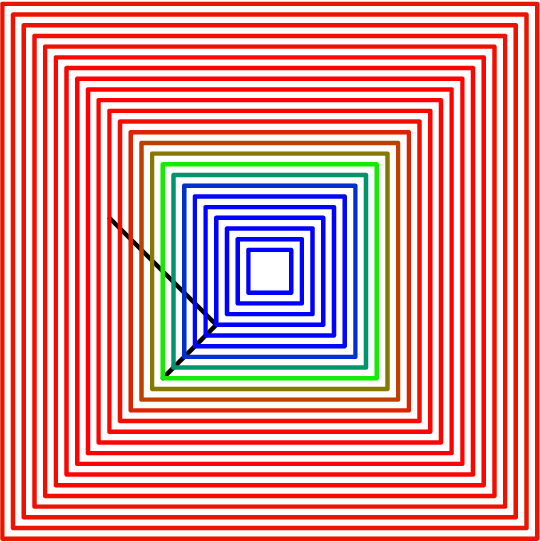}}
  \subfigure[$k_0\!=\!1,D_{P=(1,0)}^{(m)}$]{
  \includegraphics[width=0.18\textwidth]{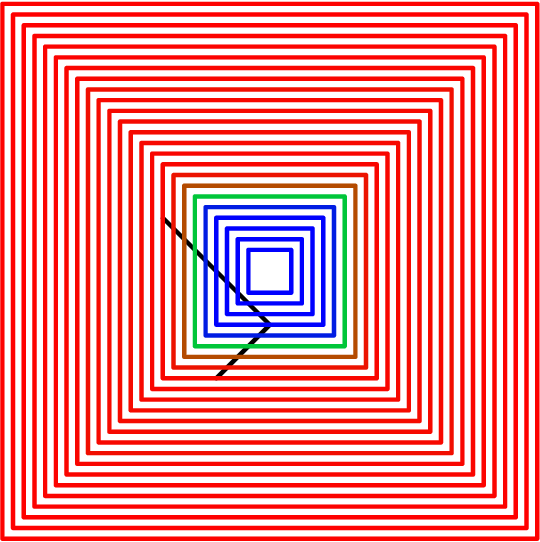}}
  \subfigure[$k_0\!=\!2,D_{P=(0,0)}^{(m)}$]{
  \includegraphics[width=0.18\textwidth]{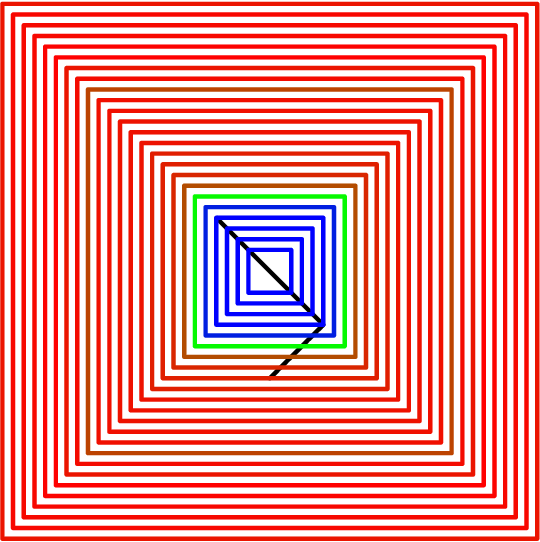}}
  \subfigure[$k_0\!=\!4,D_{P=(-1,0)}^{(m)}$]{
  \includegraphics[width=0.18\textwidth]{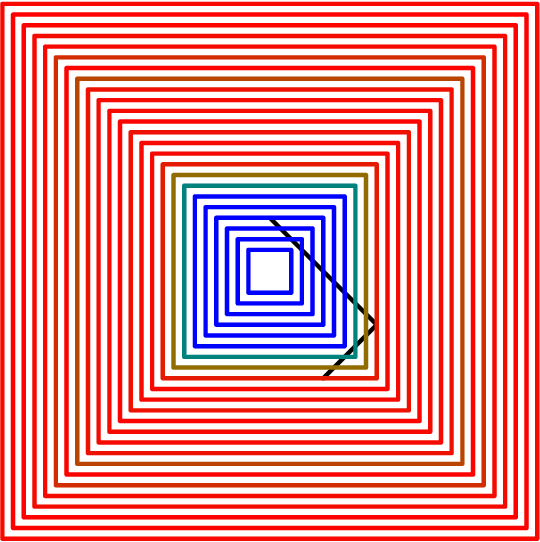}}
  \caption{Numerical results for Example \ref{example3+} {\color{xxx}with $\G$ denoted by the black solid line}.}\label{fig9}
\end{figure}
\end{example}

\begin{example}\label{ex-la}
\textbf{(Detection with limited aperture data).}
Let $\Gamma$ be a sound-soft piecewise linear crack with two tips and an interior corner given in order by $(-1,1)$, $(1,-1)$, and $(0,-2)$.
The measured limited aperture far-field {\color{xxx}data are} $\{u^\infty(\hat x_p,d_0;k_0):p=0,1,\cdots,L\}$ with $k_0=5$, $d_0=(-1,0)$, and $\hat x_p=(\cos\theta_p,\sin\theta_p)$, {\color{xxx}$\theta_p=\frac{5\pi}4+\frac{3\pi}{2L}p$, $L=800$}.

\textbf{(a) Step 1 (Detection for a rough location).}
For each $m\in\{1,\cdots,15\}$, set $D^{(m)}_P$ to be an impedance disk with impedance coefficient $\eta=ik_0$ centered at $P$ with radius $m$.
The numerical results of (\ref{211012-2}) with $D=D^{(m)}_P$ {\color{xxx}for $P\in\{10(\cos\theta,\sin\theta):\theta=-\frac\pi6,\frac\pi2,\frac{7\pi}6\}$} and $\alpha=10^{-8}$ are given by Figure \ref{fig0401}.
\begin{figure}[htb]
  \centering
  \includegraphics[width=0.3\textwidth]{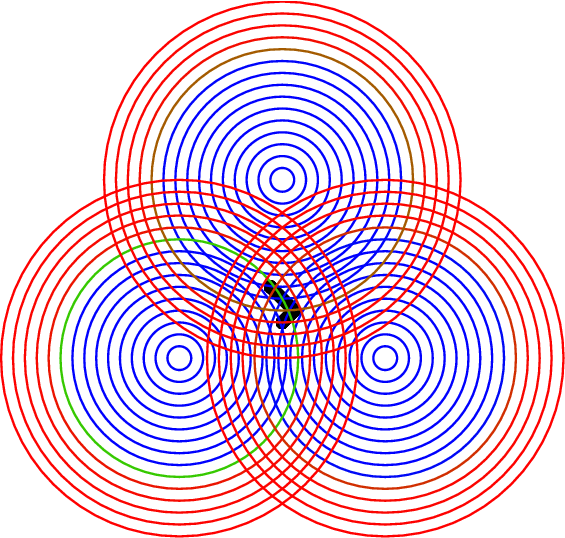}
  \caption{Numerical results for Example \ref{ex-la} (a) {\color{xxx}with $\G$ denoted by the black thick line}.}\label{fig0401}
\end{figure}

\textbf{(b) Step 2 (Detection for a precise location and the convex hull of the crack).}
Analogously to Example \ref{example3} (b), we firstly set $P=(10,0)$ and calculate the values of {\color{xxx}$\widetilde I_2(B_r(P))$ defined by} (\ref{211012-2}) with ${\color{xxx}\Sp_0}:=\{\hat x=(\cos\theta,\sin\theta):\theta\in(\frac{5\pi}4,\frac{11\pi}4)\}$ and $r=m/2$ for $m=2,3,\cdots,30$, where the regularization {\color{xxx}parameter in} (\ref{211012-2}) is {\color{xxx}set to be} $\alpha=10^{-8}$.
The numerical results for (\ref{211012-2}) are shown in Figure \ref{fig230820-1} (b)--(c) with Figure \ref{fig230820-1} (a) displaying the {\color{xxx}incident direction by the black arrow and the observation aperture by the red solid line}.
According to Figure \ref{fig230820-1} (c), the threshold is set to be $\epsilon=0.25$.
Secondly, we set $P_j=10(\cos\theta_j,\sin\theta_j)$ with $\theta_j=\frac{\pi}{16}j$, $j=0,1,\cdots,31$.
Define $m_{la,j}:=\max\{m\in\Z_+:{\color{xxx}\widetilde I_2(B_{0.1m}(P_j))}<\epsilon\}$, where the regularization {\color{xxx}parameter in} (\ref{211012-2}) is also {\color{xxx}set to be} $\alpha=10^{-8}$.
Define $\tilde\chi_{P_j,la}(x)=1$ if $|x-P_j|\leq 0.1m_{la,j}$ and $\tilde\chi_{P_j,la}(x)=0$ if $|x-P_j|> 0.1m_{la,j}$ for each $j$.
Define $\hat\chi_{P_j,la}(x)=1$ if $|x-P_j|\leq 0.1(m_{la,j}+3)$ and $\hat\chi_{P_j,la}(x)=0$ if $|x-P_j|> 0.1(m_{la,j}+3)$ for each $j$.
{\color{xxx}The numerical results for $\sum_{j=0}^{31}\tilde\chi_{P_j,la}(x)$ and $\sum_{j=0}^{31}\hat\chi_{P_j,la}(x)$ are} shown in Figure \ref{fig230820-1} (d) and (e), respectively, {\color{xxx}where the colors are given by the Matlab colormap 'jet'}.
\begin{figure}[htb]
  \centering
  \subfigure[]{
  \includegraphics[width=0.1\textwidth]{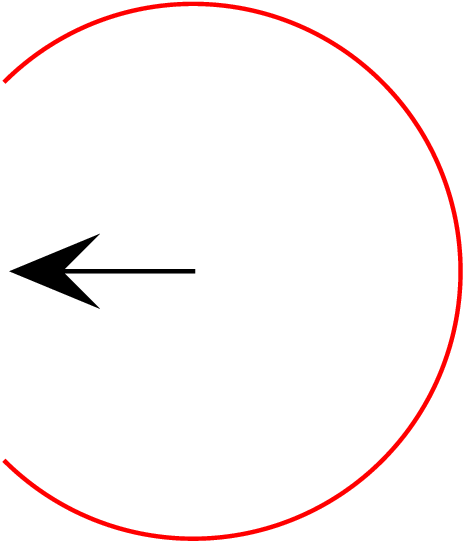}}\quad
  \subfigure[]{
  \includegraphics[width=0.18\textwidth]{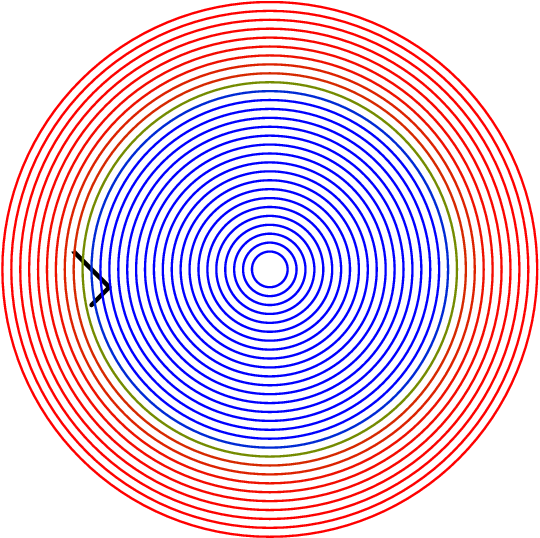}}
  \subfigure[]{
  \includegraphics[width=0.3\textwidth]{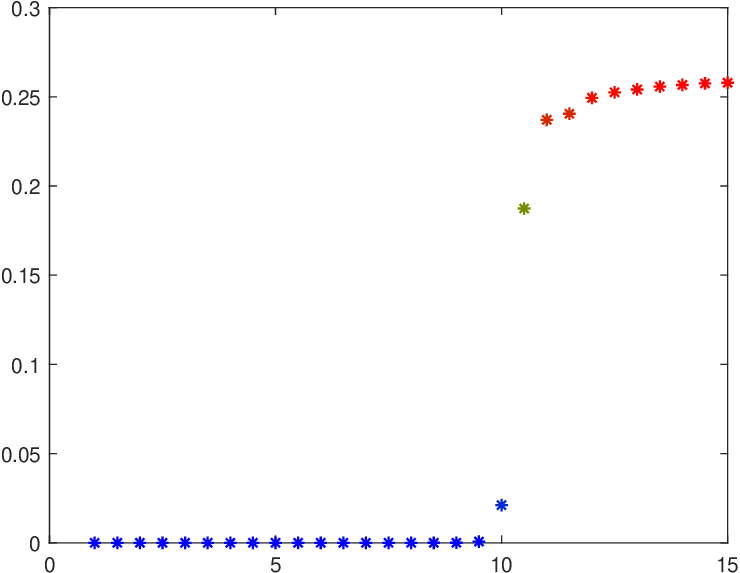}}
  \subfigure[]{
  \includegraphics[width=0.35\textwidth]{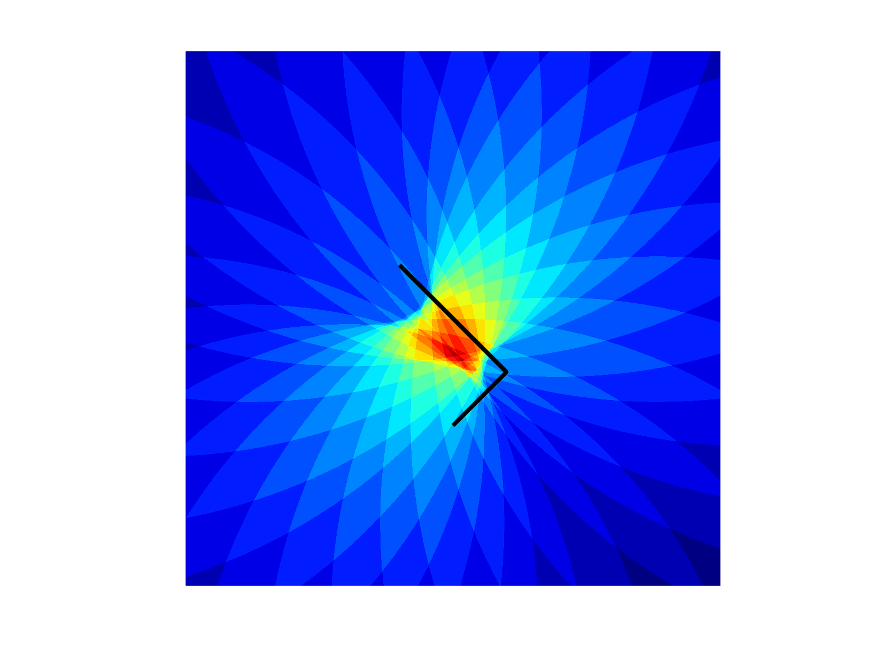}}
  \subfigure[]{
  \includegraphics[width=0.35\textwidth]{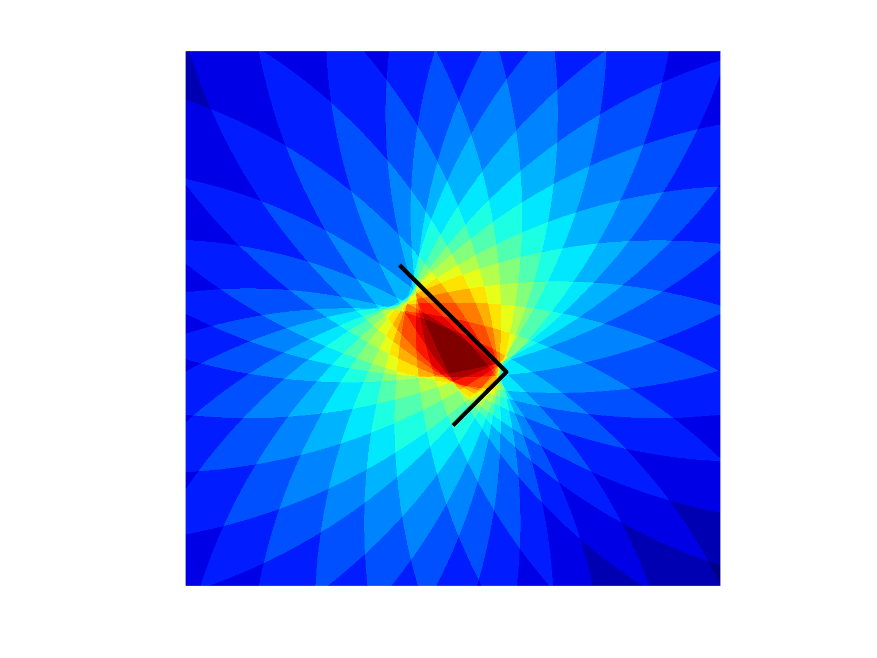}}
  \caption{Numerical results for Example \ref{ex-la} (b) {\color{xxx}with $\G$ denoted by the black solid line}.}\label{fig230820-1}
\end{figure}

\textbf{(c) Step 3 (Iteration for a precise shape).}
For a more precise numerical result, we apply the iteration method based on noisy far-field data $\{u_\delta^\infty(\hat x_p,d_0;k_0):p=0,1,\cdots,L\}$ with $k_0=5$ and $d_0=(-1,0)$, and $\hat x_p=(\cos\theta_p,\sin\theta_p)$, $\theta_p=\frac{5\pi}4+\frac{3\pi}{2L}p$, $L=40$.
The initial guess is given by the location of corners in order by $(-0.8,1.2)$, $(0.5,-1.2)$, $(-0.2,-1.9)$.
We choose different noise ratios $\delta\!>\!0$ as shown in Figure \ref{fig230903} (a)--(c).
Noting that {\color{xxx}one cannot} take it for grant that $\G$ consists of two straight lines, we also give numerical results with initial guess given by $(-0.8,1.2)$, $(-0.2,-0.2)$, $(0.7,-1.2)$, $(-0.2,-1.9)$ as shown in Figure \ref{fig230903} (d) and by $(-0.8,1.2)$, $(-0.2,-0.2)$, $(0.8,-0.8)$, $(0.3,-1)$, $(-0.2,-1.9)$ as shown in Figure \ref{fig230903} (e), respectively.
The number of total iteration steps for each figure is $10$, and we set $\alpha=20$ in (\ref{211014-1}) and $\alpha_0=10^{-2}$ in (\ref{211013-2}).

\begin{figure}[htb]
  \centering
  \subfigure[$\delta=0\%$]{
  \includegraphics[width=0.18\textwidth]{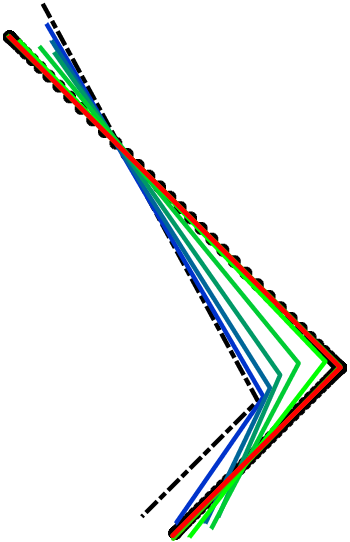}}
  \subfigure[$\delta=1\%$]{
  \includegraphics[width=0.18\textwidth]{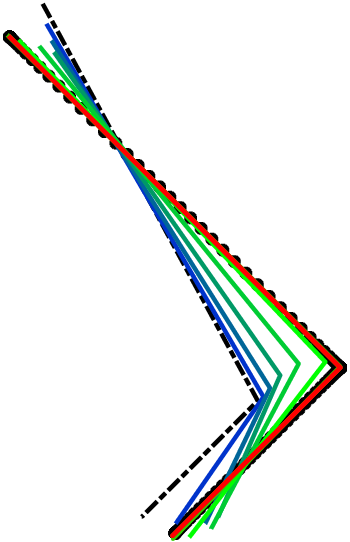}}
  \subfigure[$\delta=10\%$]{
  \includegraphics[width=0.18\textwidth]{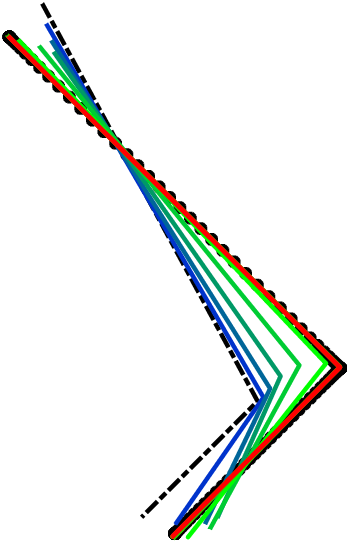}}
  \subfigure[$\delta=0\%$]{
  \includegraphics[width=0.18\textwidth]{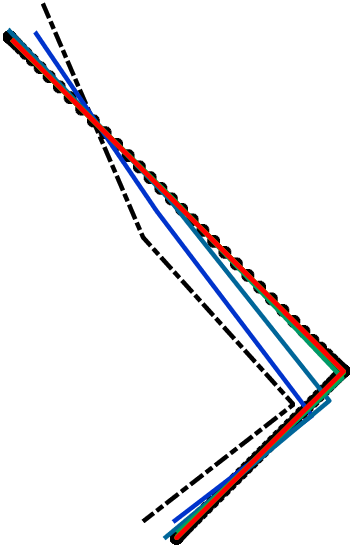}}
  \subfigure[$\delta=0\%$]{
  \includegraphics[width=0.18\textwidth]{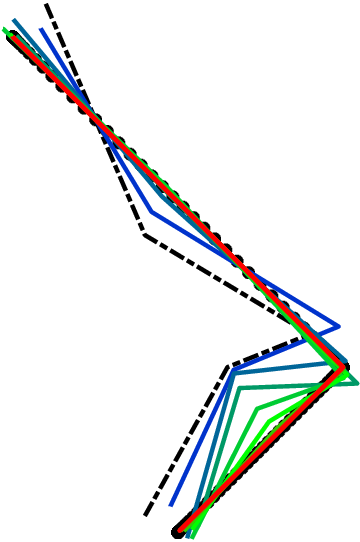}}
  \caption{{\color{xxx}Numerical results} for Example \ref{ex-la} (c). The black dots '*' represent the true crack {\color{xxx}$\G$}, the dashed line '-.' represents the initial guess, and colored lines represent the numerical result in each iteration step.}\label{fig230903}
\end{figure}
\end{example}

\section*{Acknowledgements}

The work of Xiaoxu Xu is supported by National Natural Science Foundation of China grant 12201489, the Young Talent Support Plan of Xi'an Jiaotong University, and the Fundamental
Research Funds for the Central Universities grant number xzy012022009. The work of Guanghui Hu is
supported by National Natural Science Foundation of China grant 12071236, Fundamental Research
Funds for the Central Universities in China grant 63213025.
%\bibliographystyle{siam}
%\bibliography{ref}

\end{document}